\documentclass[ijoc,nonblindrev]{informs4}

\OneAndAHalfSpacedXI

\usepackage{amsmath,amssymb,amsfonts}

\usepackage{natbib}
 \bibpunct[, ]{(}{)}{,}{a}{}{,}%
 \def\bibfont{\small}%

\usepackage{rotating}
\usepackage{fancyvrb}

\usepackage{caption}
\usepackage{booktabs} 
\usepackage{natbib}
\usepackage{subfigure}
\usepackage{graphicx} % Required for including images
\usepackage{multirow}
\usepackage{algorithm}
\usepackage{algorithmic}
\usepackage[table]{xcolor}
\usepackage{tabularx}
 \def\TheoremsNumberedThrough{%
 \theoremstyle{TH}%
 \newtheorem{theorem}{Theorem}
 
 \newtheorem{proposition}{Proposition}
 \newtheorem{corollary}{Corollary}

 \theoremstyle{EX}
 
 \newtheorem{example}{Example}
 
 \newtheorem{definition}{Definition}

 \newtheorem{observation}[definition]{Observation}
 }

\usepackage{tikz}
\tikzset{>=latex} % for LaTeX arrow head
\usetikzlibrary{positioning, arrows.meta, shapes.geometric}
\usetikzlibrary {matrix}
\tikzstyle{io} = [trapezium, 
trapezium stretches=true, % A later addition
trapezium left angle=70, 
trapezium right angle=110, 
minimum width=3cm, 
minimum height=1cm, text centered,
text width=3cm,
draw=black, fill=blue!30]

\tikzstyle{process} = [rectangle, 
minimum width=3cm, 
minimum height=1cm, 
text centered, 
text width=3cm, 
draw=black, 
fill=orange!30]
\colorlet{mylightred}{red!50!black!15}
\colorlet{mylightblue}{blue!50!black!15}
\colorlet{mylightgreen}{green!50!black!15}
\colorlet{mylightgrey}{black!20!white}

\TheoremsNumberedThrough     % Preferred (Theorem 1, Lemma 1, Theorem 2)
\ECRepeatTheorems
\JOURNAL{INFORMS Journal on Computing}

\EquationsNumberedThrough    % Default: (1), (2), ...
\MANUSCRIPTNO{}

\begin{document}
%\setcounter{page}{1}
%\DOI{10.1287/xxxx.0000.0000}%


\RUNTITLE{Learn to formulate surrogate model}

\TITLE{Learn to formulate: A surrogate model framework for generalized assignment problem with routing constraints}

\ARTICLEAUTHORS{%
\AUTHOR{Sen Xue, Chuanhou Gao}
%,\textsuperscript{a} Second Author,\textsuperscript{b} Third Author,\textsuperscript{c} Fourth Author,\textsuperscript{c}
\AFF{School of Mathematical Sciences, Zhejiang University, Hangzhou, China, \EMAIL{senx@zju.edu.cn, gaochou@zju.edu.cn}}
%\textsuperscript{b}School of Industrial Engineering, Good College, Collegeville, Maine 01234 \EMAIL{secauth@goodcoll.edu}; 
%\textsuperscript{c}Their Common Affiliation \EMAIL{thauth@anywhere.edu, fourauth@anywhere.edu}

%mirko.janc@informs.org
% \AUTHOR{Second Author}

% \AFF{School of Industrial Engineering, Good College, Collegeville, Maine 01234, \EMAIL{secauth@goodcoll.edu}}

% \AUTHOR{Third Author, Fourth Author}

% \AFF{Their Common Affiliation \{thauth@anywhere.edu, fourauth@anywhere.edu\}}
}

\ABSTRACT{%
The generalized assignment problem with routing constraints, e.g. the vehicle routing problem, has essential practical relevance. This paper focuses on addressing the complexities of the problem by learning a surrogate model with reduced variables and reconstructed constraints. A surrogate model framework is presented with a class of surrogate models and a learning method to acquire parameters. The paper further provides theoretical results regarding the representational power and statistical properties to explore the effectiveness of this framework. Numerical experiments based on two practical problem classes demonstrate the accuracy and efficiency of the framework. The resulting surrogate models perform comparably to or surpass the state-of-the-art heuristics on average. Our findings provide empirical evidence for the effectiveness of utilizing size-reduced and reconstructed surrogate models in producing high-quality solutions.}


\maketitle

\section{Introduction}

The generalized assignment problem (GAP) is a fundamental model in combinatorial optimization. It requires the determination of an optimal allocation of tasks to agents while adhering to knapsack constraints. This problem bears significant practical relevance across a broad range of industrial contexts \citep{pentico2007assignment}. For instance, in the logistics field, a basic concern is to properly allocate a set of orders to distinct vehicles for shipping. Additional applications have been developed in foundational studies on scheduling \citep{aringhieri2015two, cacchiani2017optimal, boccia2023parallel}, production planning \citep{benjaafar2004demand, dobson2001batch}, and location problems \citep{fischetti2016benders, mikic2019less}.

Many applications are formulated based on the GAP with additional routing constraints for practical needs. In this study, we specifically define these GAP-based problems as the generalized assignment problem with routing constraints (GAPR), and detail this definition in subsequent sections. A common example is the capacitated vehicle routing problem (CVRP) and its variants. This integration necessitates the inclusion of further constraints and variables to address practical considerations, such as path distance or delivery time windows. A similar scenario is observed in warehouse management, where orders are grouped into batches and picked from their respective storage locations. According to \cite{bartholdi2014warehouse}, these operations account for 55\% of the overall warehouse expenses. Recent studies increasingly investigate the inclusion of sequential picker routing in joint modeling \citep{aerts2021joint,zhang2023improved}. This integration reportedly yields more cost reduction in practical outcomes than focusing on single processes \citep{won2005joint}. Other integrated problems can be found in studies by \cite{matusiak2017utilizing, yadav2022integrated}. Consequently, these extended problems demand more complex models with increasing challenges for solutions. It further raises the challenge of solving efficiency, which is a vital measure of practicality.

Reviews on related problems show that manually designed and customized heuristics are the most widely used for the GAPR, with exact approaches being the minority \citep{pardo2024order,gutierrez2022vrp}. Exact approaches are superior in small-scale problems but are limited by the exponential growth of model size. Moreover, machine learning (ML) methods are widely integrated into optimization methods for effective problem-solving, such as branching strategy \citep{khalil2016learning,gasse2019exact}, cutting plane selection \citep{paulus2022learning}, and large neighborhood search \citep{wu2021learning}. However, little research considers the modification of the formulation structure, which might be more fundamental. This relatively leaves a gap for our study.

In this paper, we introduce a general-purpose learning framework to directly approximate the GAPR formulations. Considering the extensive diversity of routing constraints depending on reality needs, this approach highlights its advantages by automatically learning and solving without requiring problem-specific algorithm designation. Our purpose differs from the traditional exact approach that seeks a refined and tightened formulation while maintaining all necessary constraints. This approach involves the reduction of variables and the reconstruction of constraints. Supporting evidence and motivation are found in studies on branching strategies \citep{khalil2016learning,gasse2019exact}, where fewer branching nodes usually indicate less search time within the branch-and-bound (B\&B) tree. Therefore, by implementing a model with fewer integer variables, the B\&B tree might be significantly reduced. Consequently, the resultant closely-approximated model of the GAPR can be expected to produce near-optimal solutions more efficiently.

\subsection{Literature review}
The topic of learning to optimize (L2O) has attracted increasing attention in recent years. Our discussion focuses on the works on L2O related to mixed integer programming (MIP) optimization. Notably, studies by \cite{zhang2023survey}, \cite{karimi2022machine} and \cite{bengio2021machine} have investigated comprehensive reviews. A common purpose is to accelerate the MIP problem-solving process. To fulfill this purpose, parameterized ML models are trained to replace either parts or the entire computationally expensive algorithms. A representative application is the designate of the branching strategy within the B\&B algorithm. Studies by \cite{khalil2016learning} and \cite{gasse2019exact} utilized ML models to imitate the strong branching strategy, which is an exact one-step look-ahead strategy, with a high computational cost. This approximation allows a close high performance while significantly reducing the computational time required. 

However, relatively few studies concern the formulation of MIP. An aspect highlighted within these studies is that certain constraints or the objective function are not explicit. To address this issue, ML techniques are employed to approximate these undefined elements. Subsequently, the approximation ML models are embedded or integrated into the MIP formulation for further optimization. The resultant MIP models can also be referred to as surrogate models in response to the true models that cannot be explicitly expressed. Notably, \cite{schweidtmann2022optimization} and \cite{fajemisin2023optimization} have provided comprehensive overviews on the topic as \textit{constraint learning} but with a wider range of continuous and discrete optimization problems. We discuss the MIP-related works based on their ML methods.

The regression or linear approximation functions are widely used to construct MIP formulations. This method is extensively utilized in engineering with ﬁrst- and second-order polynomials. Take the work by \cite{bertsimas2016analytics} as an example. This modeling approach is widely recognized for its directness and efficiency. \cite{kleijnen2018design} presented an overview.

The tree ensembles are popular tools for constraint learning. Research conducted by \cite{biggs2017optimizing} and \cite{mivsic2020optimization} has introduced general frameworks that enable the representation of tree ensembles within MIP formulations. The applications of tree ensembles for optimization include the investigation of the electricity network problem \citep{gutina2023optimization} and the optimal power flow problem \citep{halilbavsic2018data}. Tree ensembles are valued for their interpretability for analysis and practice purposes \citep{fajemisin2023optimization}. However, the complexity of representing decision trees is a function of the tree size, making the resolution challenging with the growth of the tree.  \cite{mivsic2020optimization} proposed a Benders decomposition method for efficient problem-solving. Concurrently, the study by \cite{mistry2021mixed} provided a branch-and-cut (B\&C) algorithm with specific branching rules.

The utilization of neural networks is another widespread approach for generating learned models. The study by \cite{fischetti2018deep} initially proposed that both deep neural networks and convolutional neural networks could be presented as MIP formulations. Recent applications include the optimization of chemical processes \citep{zhou2024accelerating}, the optimal power flow problem \citep{kilwein2023optimization}, and the task scheduling problem \citep{rigo2023mppt}. The embedding of neural networks is also challenged by the large size and solving difficulty of the resultant MIP formulation \citep{fajemisin2023optimization}. To address this issue, \cite{grimstad2019relu} introduced bound-tightening algorithms to refine the values within the big-M constraints. \cite{anderson2020strong} proposed valid inequalities to strengthen the MIP formulations. Further, \cite{wang2023optimizing} utilized a Lagrangian relaxation-based B\&C procedure to produce higher-quality solutions.  

% Few studies directly use MIP formulations to develop surrogate models. The study by \cite{pawlak2017automatic} proposed a novel method to synthesize MIP constraints from feasible and infeasible solutions. Similarly, The work by \cite{kumar2019acquiring} introduced a tensor-based approach to automatically model MIP formulation for complex problems. Nonetheless, the overview by \cite{fajemisin2023optimization} highlighted a potential drawback of this approach, specifically its lack of flexibility and limited learning capacity.

Our study could provide a new perspective for L2O research. Firstly, most L2O studies, which focus on acceleration, rely on pre-established models with large feasible solution space. The resultant improvements might be therefore inherently limited. In contrast, our approach modifies and simplifies the formulations, leading to reduced search space and potential further improvements in efficiency. Secondly, many ML-based methods encounter obstacles in solving speed, necessitating further optimization efforts. A possible reason might be that the ML techniques are initially built for general learning tasks, without adaptations for the optimization of their embedded models. Our framework is specifically designed with elements of GAP, allowing for optimization without additional modifications. To the best of our knowledge, no study has been conducted on learning-based formulation modification on generalized assignment problems.

\subsection{Contributions and structure}
Our contributions in this paper are listed as follows:

1. We present a direct approximation framework with a class of surrogate models and a learning method for the efficient problem-solving of the GAPR.  

2. We provide theoretical results on the effectiveness of our framework on aspects of learning capacity and statistical analysis.

3. We demonstrate the efficacy of the framework by conducting numerical experiments on two classes of the GAPR of practical relevance. 

This paper is structured as follows. Section \ref{sec: Pd} presents a description and general formulations of the GAPR. Section \ref{sec: method} introduces our methodology to maintain a size-reduced and reconstructed surrogate model. Section \ref{sec: thm} provides theoretical results of the representational power and statistical properties of our framework. Section \ref{sec: exp} examines our approach comparing the exact algorithms. Finally, Section \ref{sec: con} summarizes the study and outlines future research directions.

\section{Problem description}
\label{sec: Pd}
We consider the GAPR as a problem of finding the optimal task assignment and route selection that minimizes both the allocation and transportation costs. Denote a set of tasks as $I$ and a set of agents as $J$. Consider a graph $G$, where one of its nodes serves as the depot. Each task corresponds to one or more nodes in $G$. Once a task is assigned to an agent, the agent must visit its corresponding nodes in $G$ from the depot with routing constraints, such as subroute elimination and soft time windows. 

Following the description of the GAPR, we provide a general MIP formulation. Suppose that task $i \in I$ weights $w_i$. The capacity of an agent $j$ is $Q_j$. Denote $y_{ij}$ as a binary variable for determining the assignment of task $i \in I$ to agent $j \in J$, where $y_{ij}=1$ indicates that task $i$ is assigned to agent $j$, and $y_{ij}=0$ otherwise. Let $y= (y_{11}, y_{12},\dots,y_{ij},\dots,y_{|I| |J|})$. Due to the absence of predefined routing rules, we denote the remaining variables as $u \in \mathbb{R}^m \times \mathbb{Z}^n$ with domain of $D_u$ and constraints coefficients as matrices $B \in \mathbb{R}^{r \times |I||J|}$, $E \in \mathbb{R}^{r \times m \times n}$ and $e \in \mathbb{R}^r$. Here $m$, $n$, and $r$ are arbitrary positive integers. The GAPR can be formulated as follows.
\begin{subequations}
\begin{align}
    \min  \quad c^ \top  y + & d ^\top u \label{obj} \\
    \text{s.t.}  \quad \sum_{j \in J} y_{ij} &=1, \quad \forall i \in I;\label{eq: set-partition}\\
    \sum_{i \in I} w_i y_{ij} &\le Q_{j}, \quad \forall j \in J; \label{eq: cap-y}\\
    \sum_{i \in I} y_{ij} &\ge 1, \quad \forall j \in J^*; \label{eq: nonempty} \\
    B y + E u &\le e, \label{eq: aux-y-u}\\
    y_{ij} &\in \{0,1\}, \quad \forall i \in I, \forall j \in J; \label{eq: setrangey}\\
    u &\in D_{u} \cap (\mathbb{R}^m \times \mathbb{Z}^n); \label{eq: setrangeu}
\end{align}
\label{GAPRmodel}
\end{subequations}
Here the objective function is computed as the assignment cost $c^{\top} y$ and routing cost $d ^{\top} u $. Constraints (\ref{eq: set-partition}) ensure that each task is uniquely assigned to an agent. Constraints (\ref{eq: cap-y}) guarantee that the assignment does not exceed the capacity of the agents. In certain scenarios where no agent should remain unassigned, this condition is expressed by constraints (\ref{eq: nonempty}) with the notation $J^* = J$. In cases where such a requirement is not applicable, we set $J = \emptyset$. All routing constraints are denoted by constraints (\ref{eq: aux-y-u}). Finally, the variable domains are defined by constraints (\ref{eq: setrangey}) and (\ref{eq: setrangeu}).

We further introduce two assumptions for the GAPR. Firstly, agents are assumed to be homogeneous, implying that capacity limits $Q_j$ are uniformly denoted by $Q$ across all agents. Secondly, let $y^*$ be an integer solution, representing an assignment plan, that satisfies constraints (\ref{eq: set-partition}), (\ref{eq: cap-y}), (\ref{eq: nonempty}) and (\ref{eq: setrangey}). We assume then there exists a feasible $u^*$ corresponding to $y^*$ that satisfies the constraints (\ref{eq: aux-y-u}) and (\ref{eq: setrangeu}). The assumption proposes that any assignment being feasible to the GAP constraints is at least feasible to the routing constraints. This assumption still covers most CVRP variants and order batching problems. In certain scenarios such as hard time windows, assignments could be infeasible due to routing constraints. However, early-stage studies have indicated that the hard time windows can be reformulated into soft constraints by incorporating appropriate penalties on these infeasible solutions \citep{taillard1997tabu}.

The objective value of the GAPR model can then be represented as a function of $y$ with the assumptions above, providing a learning target. Initially, we define the feasible region of $y$ as $P$ under the assignment and knapsack constraints, thereby also feasible to the GAPR as assumed. 
\begin{equation}
    P = \{ y : (\ref{eq: set-partition}), (\ref{eq: cap-y}), (\ref{eq: nonempty}), (\ref{eq: setrangey})\}.
    \label{eq: P}
\end{equation}
Therefore, the objective value can be presented as a function to minimize the remaining variables $u$ with given $\hat{y}$ as follows:
\begin{equation}
    f_{\text{obj}}(\hat{y}) = \min_{u} \{c^{\top} y + \lambda ^{\top} u: y=\hat{y}, (\ref{eq: set-partition}), (\ref{eq: cap-y}),(\ref{eq: nonempty}),(\ref{eq: aux-y-u}), (\ref{eq: setrangey}), (\ref{eq: setrangeu})\}.
    \label{eq: original}
\end{equation}
Finally, the task to minimize the GAPR model can be proposed as an optimization problem on $ f_{\text{obj}}(y)$:
\begin{equation}
    \min_{y\in P} f_{\text{obj}}(y) ,
\end{equation}
as $y$ denotes integer points in $P$. In subsequent sections, we study surrogate models to approximate $f_{\text{obj}}(y)$. The surrogate models are designed with significantly reduced computational demands. As a result, the models might quickly produce near-optimal solutions of the original GAPR model, thereby offering a practical approach to solving the GAPR optimization problem.

\section{Methodology}
\label{sec: method}
In this section, we introduce our methodology. Initially, we formulate a specific class of surrogate models, starting with their motivation from heuristic principles. Subsequently, we introduce the learning framework, including algorithms for sampling, regression, and iterative training.
\subsection{Surrogate model}
\subsubsection{Motivation from neighborhood search}

\begin{figure}[!t]
%\centering
\FIGURE
{\subfigure[]
{\resizebox{!}{5cm}{\includegraphics{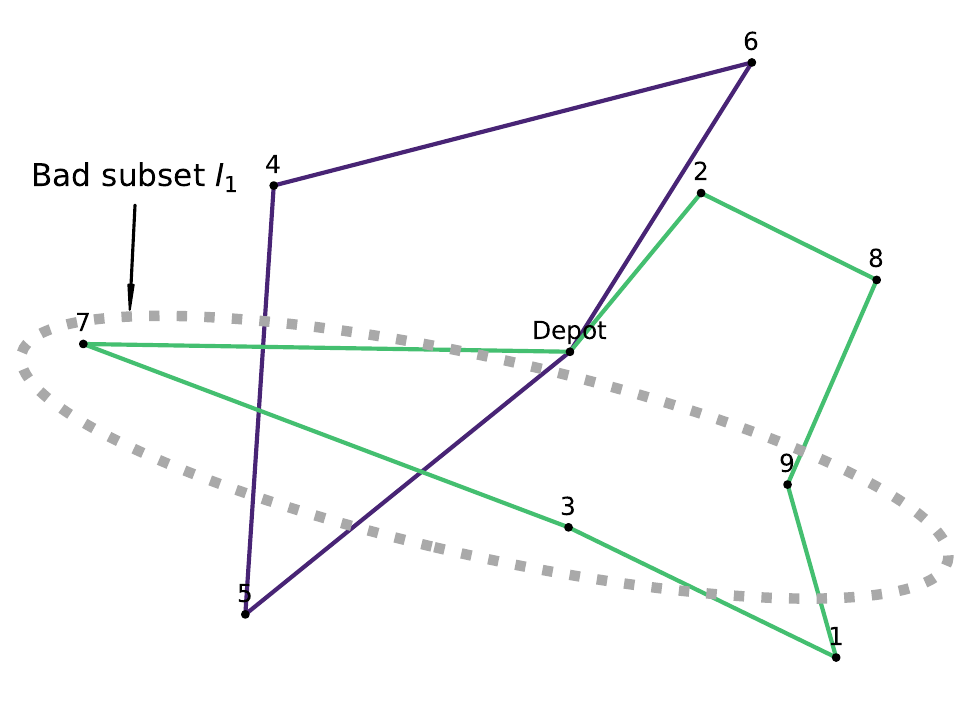}}}
\subfigure[]
{\resizebox*{!}{5cm}{\includegraphics{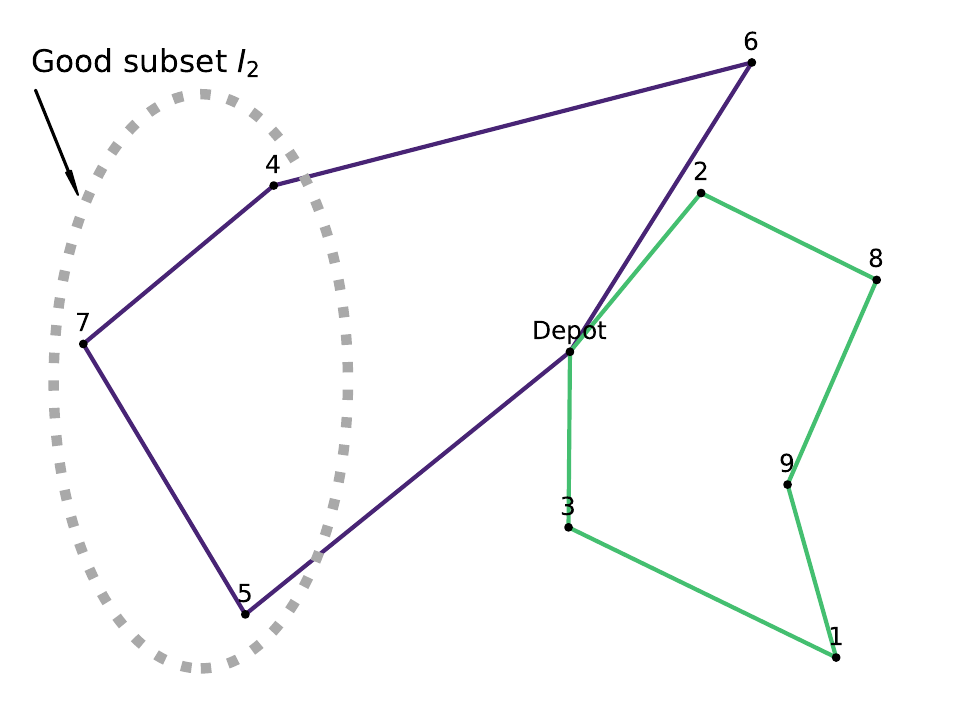}}}
}
{{A CVRP instance with two vehicles}\label{fig: instance-good-bad-subset}}
{}
\end{figure}

The development of our surrogate model is motivated by the neighborhood search (NS) strategies utilized for assignment problems. These strategies typically function by partially destroying and reconstructing a solution, while evaluating the potential improvement of such modifications. Within the context of the GAPR, we first delineate a feasible assignment as a sequence of subsets drawn from the set of tasks $I$, stated as follows:
\begin{equation}
s = (s^{k})_{k \in J}, \quad s^k \subset I \quad \forall k \in J.
\end{equation}
Here each $s^k$ represents a subset of $I$ satisfying $\bigcup_{k \in J} s^k = I$, capacity constraints, and $s^{k_1} \cap s^{k_2} = \emptyset$ for two distinct $k_1, k_2 \in J$. The NS strategies would modify parts of the solution $s$ through operations such as swapping and reallocation. To illustrate this process, we present a CVRP instance of $9$ customers with different assignment and routing paths in Figure \ref{fig: instance-good-bad-subset} as an example. In this instance, The strategies may identify subset $I_1$ in Figure (\ref{fig: instance-good-bad-subset}a) as a reason for bad performance as it initially allocates highly dispersed nodes along the same route. Conversely, The strategies may recognize $I_2$ in Figure (\ref{fig: instance-good-bad-subset}b) as an alternative good subset as it gathers closely located neighbors, indicating a lower routing distance. Given these considerations, a solution incorporating $I_1$ could be converted into one containing $I_2$ by reallocating node $7$ to be with nodes $4$ and $5$. 

However, if we could identify the good subsets in advance, the process can be transformed into a problem-solving procedure to find near-optimal solutions that contain routes including these good subsets. Denote the set of all nonempty subsets of $I$ as $\mathcal{P}^*(I)$. For example, if $I=\{1,2,3\}$, then $\mathcal{P}^*(I) = \{\{1\}, \{2\}, \{3\}, \{1,2\}, \{2,3\}, \{1,3\}, \{1,2,3\}\}$. Define $H$ as a subset of indices of $\mathcal{P}^*(I)$, where each element corresponds to a subset of $I$ denoted as $I_{\eta}$ for each $\eta \in H$. The problem of finding solutions with identified good subsets can be stated as follows.
\begin{align}
&\text{Problem: For a given set } H, \text{ find a feasible } s \text{ ensuring for any } \eta \in H, \notag \\ 
&\text{ there exists } k \in J \text{ such that } I_{\eta} \subseteq s^k.  \tag{\textbf{P1}} \label{problem1}
\end{align} 
The problem (\ref{problem1}) might produce some initially near-optimal solutions if $H$ is properly selected. To further explore this problem, we define a function $g_H$ with a given set $H$ as follows:
\begin{equation}
    g_{H}(s) = \sum_{\eta \in H} \sum_{k \in J} \textbf{1}_{I_{\eta} \cap s^k \neq \emptyset}.
\end{equation}
Here $\textbf{1}_{I_{\eta} \cap s^k \neq \emptyset}$ is a binary activation function that returns $1$ if $I_{\eta} \cap s^k$ is nonempty and $0$ otherwise. We propose the following proposition describing the relationship between $g_{H}$ and problem (\ref{problem1}).
\begin{proposition}
Let $I$ denote a set of tasks and $H$ a set of indices of $\mathcal{P}^*(I)$. Assume problem (\ref{problem1}) is feasible. Then a feasible assignment plan $\hat{s}$ is a solution to problem (\ref{problem1}) if and only if $g_{H}(\hat{s}) = |H|$.
\label{thm: prop-problem}
\end{proposition}

\proof{Proof.}
Consider two assertions to establish the necessary and sufficient conditions.

First, assume $\hat{s}$ is indeed a solution to (\ref{problem1}). Under this assumption, each subset $I_{\eta}$, for $\eta \in H$, is exclusively allocated to a specific $\hat{s}^{k(\eta)}$, denoting $k(\eta)$ is the corresponding subset in $s$ to $\eta$. This allocation implies that for any other index $k 
\neq k(\eta)$, the intersection $I_{\eta} \cap \hat{s}^k$ is empty. Consequently, the function $g_{H}$ evaluated at $\hat{s}$ can be expressed as:
\begin{equation}
    g_{H}(\hat{s}) = \sum_{\eta \in H} \sum_{k \in J} \mathbf{1}_{I_{\eta} \cap \hat{s}^k 
\neq \emptyset} = \sum_{\eta \in H} 1 = |H|.
\end{equation}

Second, assume the condition $g_{H}(\hat{s}) = |H|$ holds. Given that no $I_{\eta}$ is empty, it follows that for any $\eta \in H$, the sum $\sum_{k \in J} \mathbf{1}_{I_{\eta} \cap \hat{s}^k 
\neq \emptyset} \ge 1$. Hence, we deduce that this sum equals exactly 1 for every $\eta \in H$. This conclusion implies that each task subset $I_{\eta}$ is uniquely allocated to one subset $\hat{s}^k$. Therefore, $\hat{s}$ satisfies the exclusive allocation requirement for being a solution to (\ref{problem1}).
\endproof

The proposition \ref{thm: prop-problem} delineates an equivalence condition wherein an assignment plan $s$ is identified as a solution for the problem (\ref{problem1}). Additionally, it can be found that for any $s$, $g_{H}(s) \ge \sum_{\eta \in H} 1 = |H|$, since there is always at least one subset in $s$ contains elements in a given $I_{\eta}$.  Consequently, it enables the reformulation of (\ref{problem1}) into an optimization problem as $\min g_{H}(s)$, where $s$ denotes all feasible solutions. In the following context, we explore the property of this optimization problem in detail.

\subsubsection{Set indicator surrogate model}
We introduce an optimizable surrogate model that frames problem (\ref{problem1}) as a minimization task. Initially, we delineate the relationship between the variables $y$ in the model(\ref{GAPRmodel}) and the assignment plans $s$ as the following observation.

\begin{observation}
Let $I$ denote a set of tasks and $P$ the polyhedron of feasible integers in equation (\ref{eq: P}). The feasible assignment plans $s$ are in one-to-one correspondence to the feasible points $y$ in $P$.
\label{thm: obs}
\end{observation}
By Observation \ref{thm: obs}, we can denote the corresponding $\hat{y}$ to a specific $\hat{s}$ as $\hat{s} = s(\hat{y})$ and $\hat{y} = y(\hat{s})$ conversely. We can propose problem (\ref{problem1}) as the following optimization problem:
\begin{equation}
    \min_{y \in P} g_{H} (s(y)).
    \label{eq: min-gH}
\end{equation}
Problem (\ref{eq: min-gH}) can be further included in the following set indicator model:
\begin{subequations}
\begin{align}
    \min \sum_{\eta \in H} \beta_{\eta}& (\sum_{j \in J} \delta_{j \eta})
    \label{eq: Yobj}\\
    \text{s.t.}  \quad \sum_{j \in J} y_{ij} &=1, \quad \forall i \in I;\label{eq: Yset-partition}\\
    \sum_{i \in I} w_i y_{ij} &\le Q, \quad \forall j \in J; \label{eq: Ycap-y}\\
    \sum_{i \in I} y_{ij} &\ge 1, \quad \forall j \in J^*; \label{eq: Ynonempty} \\
    M_\eta \delta_{j \eta} &\geq \sum_{i \in I_\eta} y_{ij}, \quad \forall j \in J, \quad \forall \eta \in H; \label{eq: Ypattern}\\
    y_{ij} &\in \{0,1\},\quad \forall i \in I,\quad \forall j \in J;
    \label{eq: yrange}\\
    \delta_{j \eta} &\in \{0,1\},\quad \forall j \in J, \quad \forall \eta \in H;\label{eq: deltarange}
\end{align}
\label{SAMmodel}
\end{subequations}
Here $\beta_{\eta} \in \mathbb{R}$ and $M_\eta = |I_\eta|$ for all $\eta \in H$. Constraints (\ref{eq: Yset-partition})-(\ref{eq: Ynonempty}) are the same with Constraints (\ref{eq: set-partition})-(\ref{eq: nonempty}). The big-M constraints (\ref{eq: Ypattern}) make variables $\delta_{j \eta}$ function as a binary indicator to decide whether the allocated subset of tasks contains elements from $I_{\eta}$. 

The set indicator model is NP-hard. We can simply reduce a bin packing problem to this model. Firstly, consider an arbitrary bin packing problem with items $i \in I$ of weight $w_i$ and identical bins of capacity $Q$. The bin packing problem can be expressed as a decision problem: whether there exists a way to fit items into $|J|$ bins. This decision problem can be directly mapped onto the set indicator model with an arbitrary nonempty $H$ and $\beta_{\eta} =1$ for all $\eta \in H$. Finally, the answer to this decision problem is true if and only if the model provides a feasible solution.

The objective function of the surrogate model can also be presented as a function with $y$ being input. The surrogate model is always feasible for any $y \in P$. Therefore, by assigning a specific value to $y$, denoted as $\hat{y}$, we represent the objective value of the surrogate model as:
\begin{equation}
    \mathcal{L}_{\beta, H}(\hat{y}) = \min \{  \sum_{\eta \in H} \beta_{\eta} (\sum_{j \in J} \delta_{j \eta}): y=\hat{y}, (\ref{eq: Yset-partition}),(\ref{eq: Ycap-y}), (\ref{eq: Ynonempty}), (\ref{eq: Ypattern}),(\ref{eq: yrange}),(\ref{eq: deltarange})\},
\label{eq: LBH}
\end{equation}
where $\beta=(\beta_{\eta})_{\eta \in H}$.

The objective function $\mathcal{L}_{\beta, H}$ can be further reformulated as a linear combination of different $g_{H}$. The property enables our learning procedure to acquire the fittest $\beta_{\eta}$ values. We propose the following proposition on this property.
\begin{theorem}
Let $I$ denote a set of tasks and $P$ its resulting polyhedron of feasible integers in (\ref{eq: P}). Let $H$ denote a set of indices of $\mathcal{P}^*(I)$ and $\mathcal{L}_{\beta, H}(y)$ a function defined by equation (\ref{eq: LBH}) with $\beta_{\eta} \ge 0$ for all $\eta \in H$. Then 
\begin{equation}
     \mathcal{L}_{\beta, H}(y) = \sum_{\eta \in H} \beta_{\eta} g_{\{\eta\}} (s(y))
\end{equation} 
holds for any $y \in P$.
\label{thm: L=g}
\end{theorem}

\proof{Proof.}
For any $y \in P$, we denote $J_{\eta}$ for $\eta \in H$ as a subset of $J$ such that $s(y)^k \cap I_{\eta} \neq \emptyset$ for any $k \in J_{\eta}$. Therefore, by the definition of $g_{H}$, we have $g_{\{I_{\eta}\}} = |J_{\eta}|$. Further, under the premise of minimizing $\sum_{\eta \in H} \beta_{\eta}(\sum_{j \in J} \delta_{j \eta})$, $\delta_{j\eta}^* = 1$ for $j \in J_{\eta}$ and 0 otherwise. Thus, $\sum_{j\in J} \delta_{j\eta}^* = |J_{\eta}|$. Finally, $\sum_{\eta \in H} \beta_{\eta} (\sum_{j \in J} \delta_{j \eta})$ = $\sum_{\eta \in H} \beta_{\eta} |J_{\eta}| = \sum_{\eta \in H} \beta_{\eta} g_{\{I_{\eta}\}} (s(y))$.
\endproof

Following theorem \ref{thm: L=g}, we propose the parameterized version of $g_{H}$ as follows:
\begin{equation}
    g_{\beta,H}(s) = \sum_{\eta \in H} \beta_{\eta} g_{\{\eta\}}(s).
    \label{eq: gBH}
\end{equation}

The parameterized function $g_{\beta, H}$ extends the problem (\ref{problem1}), incorporating a stronger learning capability by assigning different weights to subsets corresponding to indices within $H$. Moreover, under the assumption that $\beta\ge \textbf{0}$, $g_{\beta,H}(s(y))= \mathcal{L}_{\beta, H}(y)$.

\subsection{Learning method}

In this subsection, we illustrate the learning method as a process to acquire the parameters $\beta$ and $H$ of the surrogate model. We first introduce sampling and regression algorithms to acquire $\beta$ with pre-defined $H$. Subsequently, we present an iterative training procedure to find $H$ and control model size.

\subsubsection{Sampling and regression algorithms}
\label{subsubsec: sampleandreg}

We illustrate the sampling method in Algorithm \ref{ag: Monte-Carlo}. The algorithm randomly selects a feasible point within the polyhedron $P$ in equation (\ref{eq: P}). This selection is achieved by setting a random assignment cost $c$ as the objective function. Given that every integer point in $P$ is an extreme point, the algorithm guarantees that every solution has a non-zero probability of being selected. Algorithm \ref{ag: Monte-Carlo} can be executed in parallel $N$ times to construct a training dataset. Assuming we have gathered $s(y_i)$ and $f_{\text{obj}}(y_i)$ for $i=1,2,\dots, N$, we denote the collection of all assignment plans as $S_{N}$ and the corresponding objective values as $b_N$.
\begin{algorithm}[H]
	\caption{Monte-Carlo sampling}
	\begin{algorithmic}[1]  
		\REQUIRE Task set $I$, agent set $J$, weights of each task $w_i$ for $i \in I$, and agent capacity $Q$ 
        \STATE {Formulate polyhedron $P$ by (\ref{eq: P}) into a MIP solver}
            \STATE {Generate random $c_{ij} \in [0,1]$ for $i\in I$ and $j \in J$}
		\STATE {Solve $\min_{y\in P} \sum_{i\in I} \sum_{j \in J} c_{ij} y_{ij}$ by the MIP solver with optimal solution as $\hat{y}$}
            \STATE {Storage $s(\hat{y})$ and $f_{\text{obj}}(\hat{y})$}
	\end{algorithmic} 
\label{ag: Monte-Carlo}
\end{algorithm}
By theorem \ref{thm: L=g}, $\mathcal{L}_{\beta, H}$ can be presented as a weighted linear combination formulation, with different $g_{H}$ being predictors. This formulation allows linear regression methods on the collected data. Therefore, to implement a regression process with a pre-defined set $H$, we initially construct the feature matrix. Let $A(S_N, H)$ represent a matrix with $N$ rows and $|H|$ columns. Denote the $j$-th element in $H$ as $h_j$. The entry in the matrix $A(S_N, H)$ located at the $i$-th row and $j$-th column is indicated as $A(S_N, H)_{ij}$ as
\begin{equation}
    A(S_N, H)_{ij} = g_{\{h_j\}} (s(y_i)).
\end{equation}
Let $\mathbb{N}_{|\mathcal{P}^*(I)|}$ be $\{1,2,\dots , |\mathcal{P}^*(I)|\}$. As a set of indices, $H$ is a subset of $\mathbb{N}_{|\mathcal{P}^*(I)|}$. Therefore, $A(S_N, H)$ is a submatrix of  $A(S_N, \mathbb{N}_{|\mathcal{P}^*(I)|})$
for any $S_N$. Below, we present a specific example for illustration. 
\begin{example}
Consider a GAPR instance with a set of items $I = \{1,2,3\}$ and a set of agents $J = \{1,2\}$. An assignment plan $s$ can be presented as an ordered 2-partition of the set $I$. Suppose we sampled three specific assignment plans, namely $s_1$, $s_2$, and $s_3$. The matrix $A(S_N, \mathbb{N}_{|\mathcal{P}^*(I)|})$ can be delineated in Table \ref{tab: example}. Furthermore, if we set $H = \{4,5,6\}$ as the indices corresponding to $\mathcal{P}^*(I)$, then the matrix $A(S_N, H)$ can be presented as the blue-colored submatrix within Table \ref{tab: example}.
\begin{table}[h]
    \TABLE
    {{The matrix $A(S_N, \mathbb{N}_{|\mathcal{P}^*(I)|})$ for three specific assignment plans, highlighting submatrix $A(S_N, H)$ with $H = \{4,5,6\}$ in blue}
    \label{tab: example}}
    {\resizebox{0.85\linewidth}{!}{%
    \begin{tabular}{c|ccccccc}
    $s$ &\{1\} &\{2\} &\{3\} &\{1,2\} & \{1,3\} & \{2,3\} & \{1,2,3\} \\ \hline
    $s_1 = (\{1\}, \{2,3\})$& 1 & 1 & 1& \cellcolor{blue!25}2 & \cellcolor{blue!25}2 & \cellcolor{blue!25}1 & 2 \\ 
    $s_2 = (\{2\}, \{1,3\})$ & 1 & 1& 1 & \cellcolor{blue!25}2 & \cellcolor{blue!25}1 & \cellcolor{blue!25}2 & 2 \\ 
    $s_3 = (\{3\}, \{1,2\})$ & 1 & 1 & 1 & \cellcolor{blue!25}1 & \cellcolor{blue!25}2 & \cellcolor{blue!25}2 & 2 \\ \hline
    \end{tabular}}}
    {}
\end{table}
\label{ex: matrix}
\end{example}
Let $b_N$ be the vector of collected objective values. Following theorem \ref{thm: L=g}, the regression problem can be expressed as a non-negative least-square as follows:
\begin{equation}
    \min_{\beta \ge \mathbf{0}}|| A(S_N, H) \beta - b_N||_2.
\label{eq: ls-original}
\end{equation}
However, in the implementation phase, we adopt a lasso regression model as follows:
\begin{equation}
    \min_{\beta} ||A(S_N, H) \beta - b_N||_2 + \gamma ||\beta||_2
    \label{eq: ls-practical}
\end{equation}
The additional regularization term controls our model size. We also drop non-negative constraints on the regression coefficients $\beta$ for speedy computing. Moreover, the resultant negative value of certain $\beta_{\eta}$ merely results in the objective value of the related $\sum_{j \in J} \delta_{j \eta}$ into a constant, without influencing the optimization outcome of the surrogate model.

\subsubsection{Training and solving procedure}

We design the training and solving process as an iterative procedure. A primary consideration here is that the set of $H$ has an exponentially large number of possibilities, rendering it impractical to incorporate all possibilities simultaneously in the regression process. Additionally, it is crucial to ensure that the surrogate model remains of a small size, considering the efficiency of optimization.

We summarize the training procedure in Algorithm \ref{ag: train}. We utilize a greedy search strategy to explore the set $H$. The search starts from a predetermined cardinality $\pi_{\text{card}}$. Consider indices of all subsets in $\mathcal{P}^*(I)$ with the cardinality of $\pi_{\text{card}}$ as potential candidates for $H$. Subsequently, a parameter $\theta$ is introduced to control the scale of the regression model in equation(\ref{eq: ls-practical}), ensuring that the size of $H$ does not exceed $\theta N$. If $H$ is within the threshold of $\theta N$, all candidate indices are included in the current $H$. If it exceeds, only the top $\theta N$ candidate indices are considered. The process progresses by increasing the $\pi_{\text{card}}$ by one if all candidate indices have been explored.

Following the determination of $H$ in each iteration, we compute the feature matrix $A(S_N, H)$ and the regularization coefficient $\gamma$. The regression model (\ref{eq: ls-practical}) is solved as an optimization problem. To evaluate the current model's quality, we compute the R-squared value. If satisfactory, we include the obtained optimal $\beta^*$ to model (\ref{SAMmodel}), yielding an assignment solution $y^*$ by optimizing. Subsequently, the objective value corresponding to $y^*$ is computed as $z$. The next process involves assessing whether this objective value is the best so far; if so, the search continues, otherwise, the program terminates. 

\begin{algorithm}[H]
\caption{Greedy search training}
\begin{algorithmic}[1]
\REQUIRE $\theta$, $\pi_{\text{card}}$, $\pi_{\text{limit}}$, $\gamma(A, b)$, $S_N$, $b_N$;
\STATE {Initialize $z^* \gets \infty$, $H \gets \emptyset$, $\Bar{H} \gets \emptyset$;}
\WHILE {$\pi_{card} \le \pi_{limit}$ }
    \IF{$\Bar{H} = \emptyset$}
        \STATE {$\Bar{H} \gets \text{indices of all subsets in } \mathcal{P}^*(I) \text{ with cardinality of } \pi_{card}$;} 
        \STATE {$\pi_{card} \gets \pi_{card} + 1$;}
    \ENDIF
    \IF{$|H| + |\Bar{H}| \le \theta N$}
        \STATE {$H \gets H \cup \Bar{H}$;}
        \STATE {$\Bar{H} \gets \emptyset$;}
    \ELSE
        \STATE {$\{\Bar{h}_1, \Bar{h}_2, \dots\} \gets \Bar{H}$;}
        \STATE {$H \gets H \cup \{\Bar{h}_1, \Bar{h}_2, \dots, \Bar{h}_{\lfloor \theta N \rfloor - |H|} \}$;}
        \STATE {$\Bar{H} \gets \{\Bar{h}_{\lceil \theta N \rceil - |H|}, \dots \}$;}
    \ENDIF
    \STATE {$\gamma \gets \gamma(A(S_N, H), b_N)$;}
    \STATE {Solve regression problem (\ref{eq: ls-practical}) and obtain optimal solution $\beta^*$;} 
    \STATE {Compute the R-squared value on the regression result as $R^2$}
    \STATE {$H \gets \{\eta: \beta_{\eta}^* \neq 0\}$;}
    \STATE {$\beta^* \gets (\beta_{\eta})_{\eta \in H}$;}
    \IF {$R^2 > R_{\text{limit}}$}
    \STATE {Solve surrogate model (\ref{SAMmodel}) built by parameters $\beta^*$ and $H$. Obtain optimal solution $y^*$;}
    \STATE {$z \gets f_{\text{obj}}(y^*)$;}
    \IF {$z^* > z$}
    \STATE {$z^* \gets z$;}
    \ELSE
    \STATE \textbf{break};
    \ENDIF
    \ENDIF
\ENDWHILE
\end{algorithmic}
\label{ag: train}
\end{algorithm}

Finally, we present the workflow of our method in Figure \ref{fig: workflow}. Given a GAPR instance, we run $N$ sampling algorithms in parallel to obtain a data set $S_N$ and $b_N$. We then feed the data into Algorithm 2 for iterative training. Then we obtain the final objective function value $z^*$. Concurrently, the optimal values of decision variables $y$ and $u$ in the model (\ref{GAPRmodel}) for practical planning can also be derived by computing the objective value.

%\begin{tikzpicture}[node distance=1.5cm and 2cm, auto, >=Latex]
\begin{figure}[h]
%\centering
\FIGURE
{
\begin{tikzpicture}[node distance=1cm and 1.5cm, auto, >=Latex, scale=0.6, every node/.style={scale=0.6}]
% Nodes 
\node (input) [io] {INPUT: A GAPR problem instance};
%\node (step1) [process, right=of input] {Step 1};
\node (process2) [process, right=0.6cm of input] {Monte-Carlo sampling Algorithm \ref{ag: Monte-Carlo}};
\node (process1) [process, above =0.25 cm of process2] {Monte-Carlo sampling Algorithm \ref{ag: Monte-Carlo}};
\node (dots) [ below =0.1 cm of process2] {\vdots};
\node (process3) [process, below=0.25cm of dots] {Monte-Carlo sampling Algorithm \ref{ag: Monte-Carlo}};
\node (step3) [process, right=2cm of process2] {Collect data into $S_N$ and $b_N$};
\node (step4) [process, right=0.6cm of step3] {Greedy search training Algorithm \ref{ag: train}};
\node (output) [io, right=0.6cm of step4] {OUTPUT: Objective value $z^*$}; 
% Arrows with arcs and text
\draw[->] (input.east) to[bend left=10] node[above] {} (process1.west);
\draw[->] (input.east) to node[below] {} (process2.west);
\draw[->] (input.east) to[bend right=10] node[below] {} (process3.west);
\draw [->] (process1.east) to[bend left=15] node[above=0.3cm] {$s(\hat{y_1}), f_{\text{obj}} (\hat{y_1})$} (step3.west);
\draw[->]   (process2) -- node[above] {$s(\hat{y_2}), f_{\text{obj}} (\hat{y_2}) $} (step3.west);
\draw[->]  (process3.east) to[bend right=10] node[below=0.5cm] {$s(\hat{y_N}), f_{\text{obj}} (\hat{y_N})$} (step3.west);
%\draw[arrows = {-Stealth[inset=0pt, length=5pt, angle'=20]}]  (process3.east) to[bend right=10] node[below] {$s(\hat{y}), f_{\text{obj}} (\hat{y}) $} (step3.west);
\draw[->] (step3) -- node[above] {} (step4);
\draw[->] (step4) -- node[above] {} (output);
\end{tikzpicture}
}
{{The flow chart of learning method}\label{fig: workflow}}
{}
\end{figure}

\section{Theoretical analysis}
\label{sec: thm}
This section theoretically explores the effectiveness of our framework. We first examine its representational power, which is a property closely associated with learning capacity. We further provide statistical estimation results to validate the effectiveness of our surrogate models to produce solutions as an approximation to the original formulations
\subsection{Representational power}
Our model utilizes a finite set of features as predictors, raising questions regarding its learning capacity. The central inquiry pertains to the model's ability to accurately fit varying objective functions specific to the GAPR problems. This concept has been discussed in studies of neural networks, particularly for graph neural networks (GNN)  \citep{azizian2021expressive}. The study by \cite{chen2022representing} further extended this discussion into learning the objective values of general MIP by GNN while delineating representational power into two fundamental aspects: separation power and approximation power. Subsequently, we follow the scheme of \cite{chen2022representing} by providing theoretical results on our method.

The separation power can be initially stated as follows: an equivalence relation between two inputs exists if, and only if, any given function produces identical outputs for these inputs. For the GAPR, we first define the equivalent relationship as follows.

\begin{definition}
    Two solutions $y_1, y_2 \in P$ are equivalent if there exists a permutation $\phi$ such that $s(y_1) = \phi(s(y_2))$. 
\end{definition}
Then we show that our surrogate model has separation power by the following theorem.
\begin{theorem}
    Two solutions $y_1, y_2 \in P$ are equivalent
    if and only if $\mathcal{L}_{\beta, H}(y_1) = \mathcal{L}_{\beta, H}(y_2)$ for any $\beta$ and $H$.
    \label{thm: separate}
\end{theorem}

\proof{Proof.}
We use the condition of equivalent and its converse to prove the necessary and sufficient conditions.

(1). If $y_1$ and $y_2$ are equivalent, then for any $\beta$ and $H$, denote $\delta^{(1)}_{j \eta}$ and $\delta^{(2)}_{j \eta}$ as the optimal solution of model (\ref{SAMmodel}) with fixing $y$ to $y_1$ and $y_2$. If $\beta_{\eta}$ is non-negative, then 
\begin{align}
    \sum_{j \in J} \delta^{(1)}_{j \eta} &= g_{\{I_{\eta}\}} (s(y_1))\\
    &=\sum_{j \in J} \textbf{1}_{I_{\eta}\cap s(y_1)^j \neq \emptyset}\\
    &=g_{\{I_{\eta}\}} (\phi(s(y_1)))\\
    &=  \sum_{j \in J} \delta^{(2)}_{j \eta}.
\end{align}
If $\beta_{\eta}$ is negative, then $ \sum_{j \in J} \delta^{(1)}_{j \eta} = \sum_{j \in J} \delta^{(2)}_{j \eta}=|J|$. Finally, we obtain that  $\mathcal{L}_{\beta, H}(y_1)  = \mathcal{L}_{\beta, H}(y_2)$ for any $\beta$ and $H$.

(2). Conversely, assume $y_1$ and $y_2$ are not equivalent. we set $\beta_{\eta} = 1$ for any $\eta$. This ensures that for any feasible input $y$, $\mathcal{L}_{\beta, H}(y) = g_{H} (s(y))$, where $s(y)$ denote the corresponding assignment plan $s$ for $y$. We define $s_1 = s(y_1)$ and $s_2 = s(y_2)$ for simplicity.

We can propose that there exist an index $k_1$ such that $s_1^{k_1} 
\neq s_2^{k_2}$ for all $k_2 \in J$. To contradict this assertion is to suggest that for every $k_1 \in J$, there exists a corresponding $k_2$ such that $s_1^{k_1} = s_2^{k_2}$. Such a scenario implies that $s_1$ and $s_2$ are merely different permutations of identical subsets of $I$. However, this would imply $f(y_1) = f(y_2)$ under the assumption of homogeneity among agents, contradicting our initial condition. Given the existence of such a $k_1$, we consider two scenarios.

Firstly, if there exists a $k_2$ such that $s_1^{k_1}$ is a subset of $s_2^{k_2}$, we define $H=\{\eta: I_{\eta} = s_2 ^{k_2}\}$. In this scenario, elements of $s_2^{k_2}$ are found in more than one subset within $s_1$, leading to $\mathcal{L}_{\beta, H}(y_2) = 1< \mathcal{L}_{\beta, H}(y_1)$.

Conversely, $s_1^{k_1}$ does not form a subset of any subset in $s_2$, implying its elements are distributed across multiple subsets in $s_2$. we define $H=\{\eta : I_{\eta} = s_1 ^{k_1}\}$. Thus, $\mathcal{L}_{\beta, H}(y_1) = 1< \mathcal{L}_{\beta, H}(y_2)$.

Therefore, there exist some $\mathcal{L}_{\beta, H}(y_1)$ that output differently.
\endproof

Upon proving the separation power, we can further show that in the best scenario, our surrogate model is able to identify a solution that is either identical to or equivalent to that derived from the original model.
\begin{corollary}
Let $(y^*, u^*)$ denote the optimal solution to a given instance of GAPR defined by model \ref{GAPRmodel}. There exists a surrogate model represented by $\mathcal{L}_{\beta, H}$ such that the equation $y = \argmin_{y}\mathcal{L}_{\beta, H}(y)$ holds if and only if $y$ is equivalent to $y^*$. 
\end{corollary}
\proof{Proof.}
Let $s^* = s(y^*)$. Set $H = \{\eta \mid \exists k \in J \text{ such that } I_{\eta} = s^{*k}\}$ and $\beta = \mathbf{1}$. According to Theorem \ref{thm: L=g}, $\mathcal{L}_{\beta, H}(y)=g_{\beta, H}(s(y)) = \sum_{\eta \in H} g_{{\eta}}(s(y))$. Here, $g_{{\eta}}(s(y)) = \sum_{k \in J} \mathbf{1}_{I{\eta} \cap s(y)^k \neq \emptyset} \ge 1$, leading to $\mathcal{L}_{\beta, H} \ge |H|$. It can be computed that $\mathcal{L}_{\beta, H}(y^*)=g_{\beta, H}(s^*) = |H|$. Thus $y^*$ is the minimal solution of $\mathcal{L}_{\beta, H}$.

Further, by Theorem \ref{thm: separate}, any $y$ which is equivalent to $y^*$ also minimizes $\mathcal{L}_{\beta, H}$. Conversely, if $y$ is not equivalent to $y^*$, it implies a discrepancy of elements in sequences $s(y)$ and $s(y^*)$. This discrepancy guarantees the presence of some $\eta$ for which $g_{{\eta}}(s(y)) \ge 2$, thereby disqualifying $y$ as the minimal solution of $\mathcal{L}_{\beta, H}$. This completes the proof.
\endproof

The approximation power involves the error between the surrogate model and the learning target. We provide an upper bound on the training error of equation (\ref{eq: ls-practical}) for any pre-defined $H$ without the regulation term. This result generally applied to all GAPR problems.

\begin{theorem}
Given a set of tasks $I$ and a set of agents $J$, with an overdetermined matrix $A(S_N, H)$, and the normalized $\Bar{b}_N = \frac{b_N}{||b_N||}$. Let $A=A(S_N, \mathbb{N}_{|\mathcal{P}^*(I)|})$ and $\lambda_{max}$ denote the largest eigenvalue of $\mathbf{Id}-A(A^\top A)^{-1}A^\top$, where $\mathbf{Id}$ is the unit matrix. It follows that 
$$\min_{\beta \in \mathbb{R}^{|H|}, H\subseteq \mathcal{P}^*(I)} ||A(S_N, H) \beta - \Bar{b}_N|| \le \sqrt{\lambda_{max}}$$.
\label{thm: error-bound}
\end{theorem}
\proof{Proof.}
We begin by noting, as illustrated in Example \ref{ex: matrix}, that the matrix $A(S_N, H)$ is a submatrix of $A$. The selection of columns from $A$ to form $A(S_N, H)$ can also be controlled by the vector $\beta$ by setting certain values to be zeros. Thus, we have
\begin{equation}
\min_{\beta \in \mathbb{R}^{|H|}, H\subseteq \mathcal{P}^*(I)} \left\|A(S_N, H) \beta - \bar{b}_N\right\| = \min_{\beta \in \mathbb{R}^{|\mathcal{P}^*(I)|}} \left\| A \beta - \bar{b}_N\right\|.
\label{eq: normlambda}
\end{equation}
The optimal solution to this least squares problem is given by $\beta^* = (A^\top A)^{-1} A^\top \bar{b}_N$.

Proceeding, we find that
\begin{align*}
    \left\| A \beta^* - \bar{b}_N\right\|^2 &= \bar{b}_N^\top \bar{b}_N - 2 \bar{b}_N^\top A \beta^* + (A \beta^*)^\top (A \beta^*)\\
    &= \bar{b}_N^\top \left( \mathbf{Id} - A (A^\top A)^{-1} A^\top \right) \bar{b}_N.
\end{align*}
We denote $\mathbf{Id} - A (A^\top A)^{-1} A^\top$ by $W$. Thus, $\left\| A \beta^* - \bar{b}_N\right\|^2 = \bar{b}_N^\top W  \bar{b}_N$. Given that $\bar{b}_N^\top W  \bar{b}_N \ge 0$ for any $\bar{b}_N$, it follows that $W$ is positive semidefinite. Furthermore, as $A (A^\top A)^{-1} A^\top$ is symmetric, $W$ is also symmetric. We have $W= U \Lambda U^\top$, where $U$ is an orthogonal matrix and $\Lambda$ is a diagonal matrix with the eigenvalues of $W$ on its diagonal. For any $\bar{b}_N \in \mathbb{R}^N$ with $\left\|\bar{b}_N\right\|=1$, we have $\bar{b}_N^\top W  \bar{b}_N = \bar{b}_N^\top U \Lambda U^\top  \bar{b}_N$. Setting $z=U^\top \bar{b}_N$ yields $\left\|z\right\|=1$ and $\bar{b}_N^\top W  \bar{b}_N = z ^\top \Lambda z \le \lambda_{\max} \left\|z\right\| = \lambda_{\max}$, with $\lambda_{max}$ being the largest eigenvalue on the diagonal of $\Lambda$  Therefore, we conclude that
\[
\left\| A \beta^* - \bar{b}_N\right\| \le \sqrt{\lambda_{\max}}.
\]
Together with equation (\ref{eq: normlambda}), we have the theorem result.
\endproof

\subsection{Statistical analysis}

We further propose theoretical results for estimating the objective value that our surrogate model yields. The key idea is to use statistical tools to analyze the collected samples. Given that the surrogate model is built from data, it is challenging to derive a general approximation ratio analysis as to the traditional approximation algorithms. However, it remains possible to estimate the objective value based on statistical inference and the optimization results of the surrogate model.

The utilization of statistical tools involves a new perspective: considering the objective value of an MIP formulation as a random variable. Despite that the feasible solutions for an MIP are usually finite, their massive amount allows the investigation of statistical properties, including statistical distribution. The exploration of statistical properties of solutions to optimization problems has been discussed in many studies. For GAP, a study by \cite{vittorio2021matheuristics} noted that on a small GAP instance, the distribution representing the objective values of all possible solutions can be fitted by a binomial distribution. For the GAPR, we consider the predefined assignment plan $s$ as an event. If an event $s$ happens, the value of the random variable is given by a specific objective function $f(y(s))$, such as the function $f_{\text{obj}}$ in equation (\ref{eq: original}). The detailed definition is as follows.
\begin{definition}
Let  $(\Omega, \mathcal{F}, \mathbb{P})$ be a probability space, where the sample space $\Omega$ is the set of all feasible assignment plans $s$, the event space $\mathcal{F}$ is the power set of $\Omega$, and $\mathbb{P}(s)$ is the probability for $s$ to be sampled. Denote the objective function as $f$. Define random variable $X^{f}(s)$ as mapping
$$  s \mapsto   f(y(s)),$$
where $f$ is a real-valued function.
\label{thm: def}
\end{definition}

\begin{figure}
\FIGURE
{\subfigure[Joint order batching and picker routing problem]
{\resizebox{!}{4.5cm}{\includegraphics{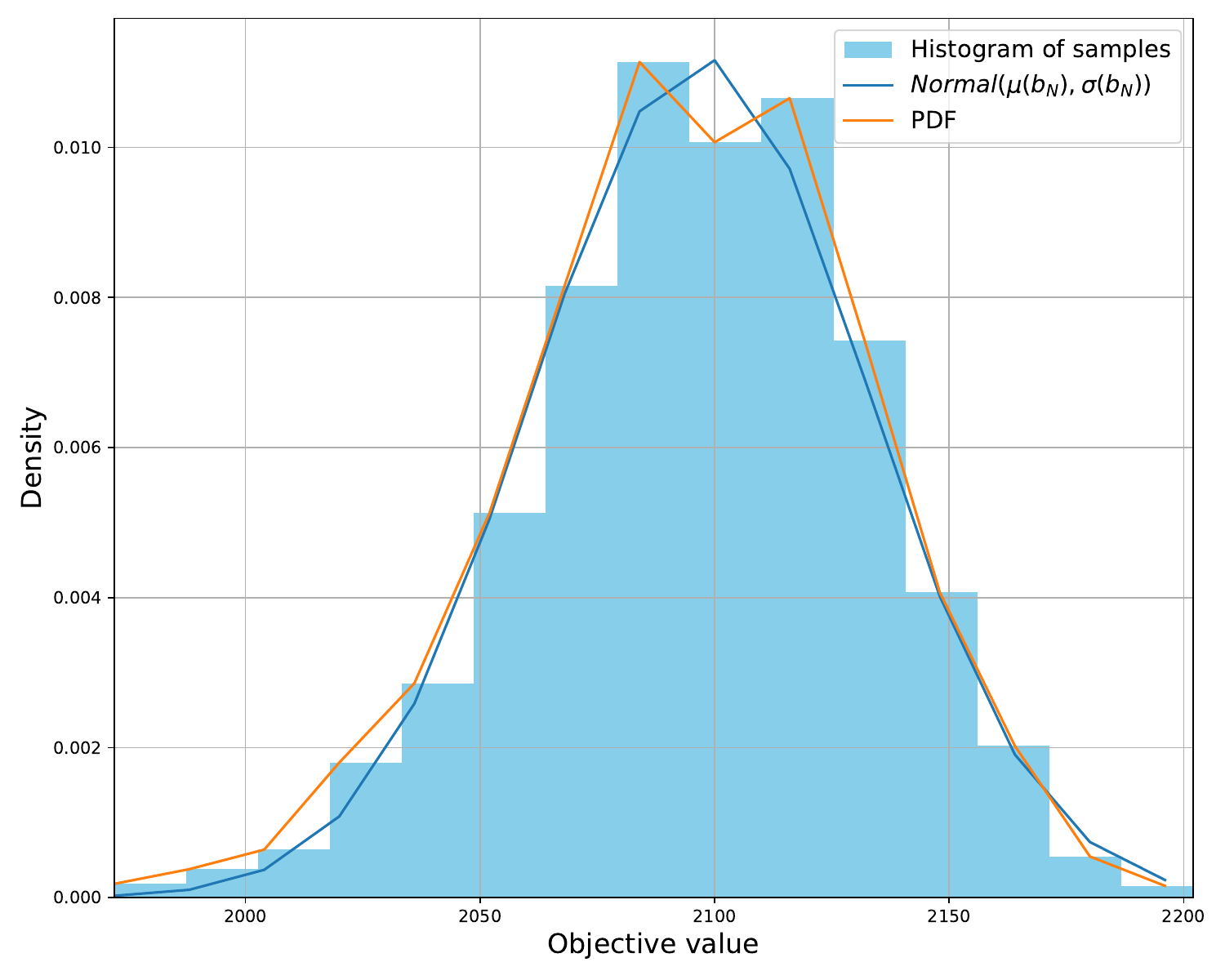}}}
\subfigure[soft-clustered vehicle routing problem]
{\resizebox*{!}{4.5cm}{\includegraphics{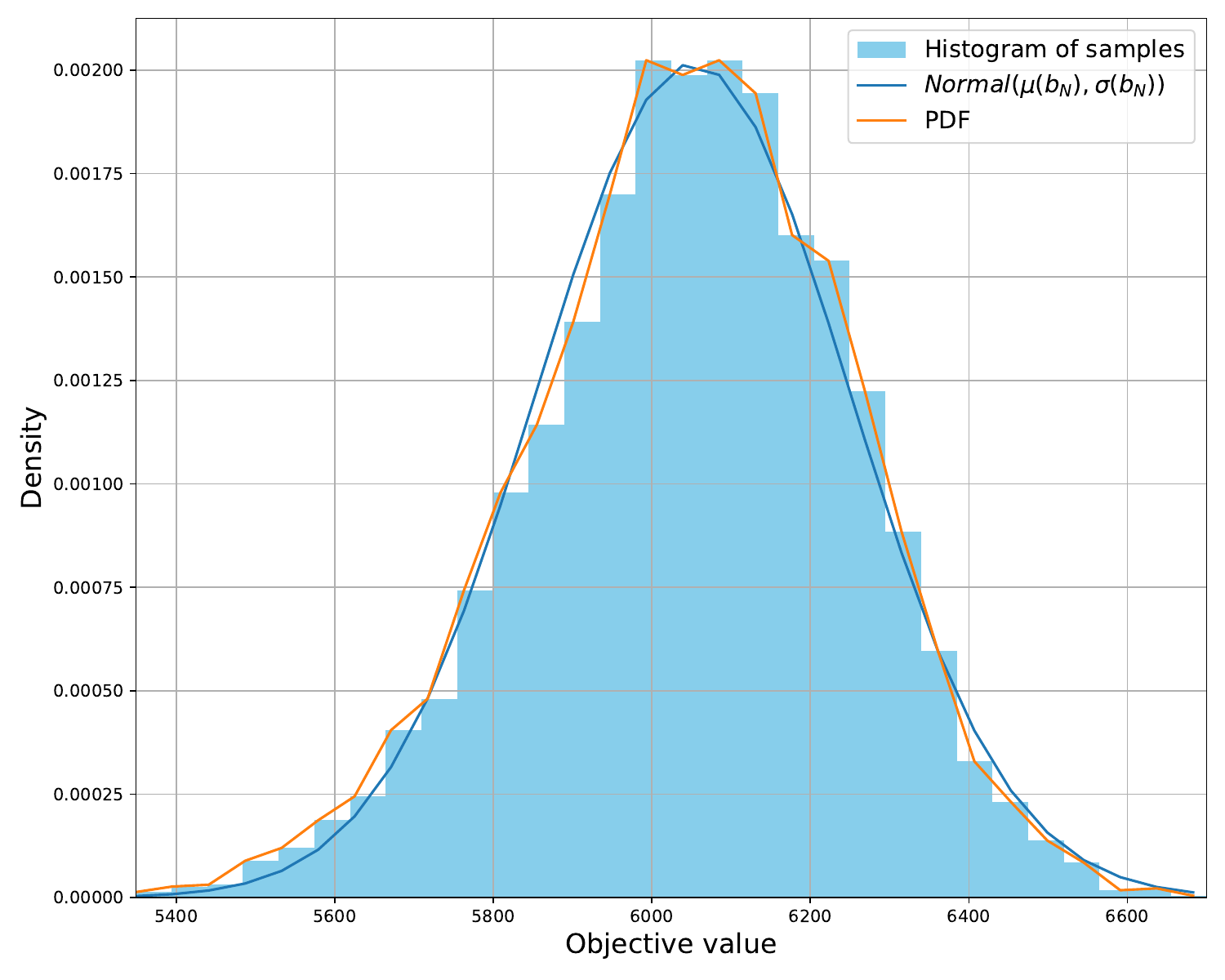}}}
}
{{Solution distribution of the GAPR problems introduced in Section \ref{sec: exp}} \label{fig: distribution}}
{}
\end{figure}
We first infer the distribution of random variables based on the objective function. To simplify notations, we denote the random variable defined by the objective function (\ref{eq: original}) as $X^f$, and that defined by the function $g_{\beta, H}$ as $Y^g$. Experiments are conducted on two distinct classes of the GAPR problems that are described in detail in Section 5. From each instance, a collection of 6,000 samples of objective values is sampled. These samples are subsequently represented in bar charts, alongside the empirical probability density function (PDF) to closely approximate their true distribution curve. Two typical plots are exhibited in Figure \ref{fig: distribution}. Moreover, we conduct Chi-squared tests to check whether or not the PDFs follow the normal distribution. The results all indicate that the samples can be fitted by normal distributions.

The distribution tests provide strong support to the assumption the distributions of $X^f$ and $Y^g$ can be approximated by the normal distribution for subsequent estimations. Under this premise, We consider the scenario upon training a surrogate model in the form of model (\ref{SAMmodel}) and solving it to a solution denoted as $\hat{y}$, yielding $\hat{Y} = g_{\beta, H}(s(\hat{y}))$. Therefore, We can estimate the expected value of the true value of $f_{\text{obj}}(\hat{y})$, without computing an optimization problem. We consider the subsequent theorem.

\begin{theorem}
\label{thm: bivariate}
Let $X^f$ and $Y^g$ be random variables defined by the definition \ref{thm: def}. If $(X^f,Y^g)$ follow bivariate normal distribution, then
\begin{equation}
\centering
(X^f,Y^g) \sim \text{Normal}\left(\left(\mu_{X}, \mu_{Y}\right),\left[\begin{array}{cc}
	\sigma_{X}^{2} & \rho \sigma_{X} \sigma_{Y} \\
	\rho \sigma_{X} \sigma_{Y} & \sigma_{Y}^{2}
\end{array}\right]\right)
\label{jointbinormal}
.\end{equation}
Given $Y^g=\hat{Y}$, the conditional expectation of $X^f$ is
\begin{equation}
    \mathbb{E}[X^f \mid Y^g=\hat{Y}]=\mu_{X}+\sigma_{X} \rho \frac{\hat{Y}-\mu_{Y}}{\sigma_{Y}}. 
    \label{eq: expectation}
\end{equation}
The upper bound $\zeta$ that $P(X^f < \zeta| Y^g=\hat{Y})=1-\alpha$ is 
\begin{equation}
 \zeta = \sigma_{X|Y} \Phi ^{-1}(1-\alpha)+\mu_{X|Y},
 \label{eq: upperboundexpectation}
\end{equation}
where $ \mu_{X|Y}=\mu_{X}+\sigma_{X} \rho \frac{\hat{Y}-\mu_{Y}}{\sigma_{Y}}$, $ \sigma_{X|Y}=(1-\rho^{2}) \sigma_{X}^{2}$ and $\Phi(\zeta)=\frac{1}{\sqrt{2 \pi}} \int_{-\infty}^{\zeta} e^{-t / 2} d t$.
\end{theorem}
\proof{Proof.}
The statistical properties utilized in this theorem can be found in multivariate statistical analysis \citep{hardle2007applied}. By the assumption of bivariate normal distribution, the conditional probability density function of $X^f$ is
$$
f_{X \mid Y}\left(X^f \right)=\frac{f_{Y,X}\left(Y^g, X^f\right)}{f_{Y}\left(Y^g\right)},
$$
where
\begin{align*}
	&f_{Y, X}\left(Y^g, X^f\right)=\\
	&\frac{1}{2 \pi \sigma_{Y} \sigma_{X} \sqrt{1-\rho}} \exp \left(-\frac{\left(\frac{\left(Y^g-\mu_{Y}\right)^{2}}{\sigma_{Y} ^{2}}+\frac{\left(X^f-\mu_{X}\right)^{2}}{\sigma_{X}^{2}}-\frac{2 \rho\left(Y^g-\mu_{Y}\right)\left(X^f-\mu_{X}\right)}{\sigma_{Y} \sigma_{X}}\right)}{2\left(1-\rho^{2}\right)}\right),  
\end{align*}
and
$$
f_{Y}(Y^g)=\frac{1}{\sqrt{2 \pi \sigma_{Y}^{2}}} \exp \left(-\frac{\left(Y^g-\mu_{Y}\right)^{2}}{2 \sigma_{Y}^{2}}\right).
$$
Therefore, the conditional distribution follows
$$
X^f \mid Y^g=\hat{Y} \sim Normal\left(\mu_{X}+\rho \frac{\sigma_{Y}}{\sigma_{X}}\left(\hat{Y}-\mu_{Y}\right),\left(1-\rho^{2}\right) \sigma_{X}^{2}\right)
.$$

This distribution directly shows the expected value of $X^f$. Moreover, based on the upper bound as
$$P(X^f < \zeta| Y^g=\hat{Y})=1-\alpha,$$ we have
$$ \Phi(\frac{\zeta-\mu_{X|Y}}{\sigma_{X|Y}})=1-\alpha.$$
Then, with the known upper critical value we have the result
$$ \zeta = \sigma_{X|Y} \Phi ^{-1}(1-\alpha)+\mu_{X|Y}.$$
\endproof

Assume our regression process yields a statistically significant correlation coefficient $\rho$. The obtained value $\hat{Y}$, as a minimized objective value in the surrogate model (\ref{SAMmodel}), diverges from its average $\mu_{Y}$. Consequently, by equation (\ref{eq: expectation})-(\ref{eq: upperboundexpectation}), the value of $X^f$ can be expected to diverge downwards from its average $\mu_{X}$, indicating a potentially high-quality solution.

% The following theorem is left to the situation that $(X^g, Y^g)$ does not follow a bivariate normal distribution. Therefore, we estimate the upper bound of the confidence interval of $X^f$ directly based on the empirical distribution function.

% \begin{theorem}
% Let $X^f$ and $Y^g$ be random variables defined by the definition \ref{thm: def}. Consider the sequence of independent and identically distributed (i.i.d.) pairs $\{(X_{i}, Y_{i})\}_{i=1}^{N}$, each following the joint cumulative distribution function $G(u, v)$. The joint empirical cumulative distribution function based on these observations is given by
% $\widehat{G}_{N}(u, v) = \frac{1}{N} \sum_{i=1}^{N} \mathbf{1}_{\{X_{i} \leq u, Y_{i} \leq v\}}$. The distribution of $Y$ is given by $\widehat{G}_{N}(v) = \frac{1}{N} \sum_{i=1}^{N} \mathbf{1}_{\{Y_{i} \leq v\}}$. Here $\mathbf{1}_{(\cdot)}$ denotes the indicator function. Consider $Y^g \leq \hat{Y}$ and a real number $\zeta$. Then

% $$ P(X^f \leq \zeta | Y^g \leq \hat{Y})= \frac{\widehat{G}_{N}(\zeta,\hat{Y})}{\widehat{G}_{N}(\hat{Y})} =\frac{\sum_{i=1}^{N} \mathbf{1}_{(X_{i} \leq \zeta, Y_{i}\leq \hat{Y})}}{\sum_{i=1}^{N} \mathbf{1}_{Y_{i}\leq \hat{Y}}}.$$
% \end{theorem}

\section{Numerical experiments}
\label{sec: exp}
In this section, we examine our methodology through numerical experiments conducted on two distinct GAPR problems. We address the following experimental questions. \textbf{Q1}) How accurately do our proposed surrogate models approximate solutions when compared to exact models? \textbf{Q2}) Do our surrogate models significantly enhance solving efficiency? Given that our surrogate models are NP-hard, it is important to empirically verify the improvement in computational efficiency. \textbf{Q3}) Can the observed improvements be attributed to the modifications introduced in the MIP formulations? This question investigates the causality behind improvements to further enhance the understanding of the contribution of our methodology. The notations for data presentation are listed in Table \ref{tab: notations}.

Finally, to further demonstrate the practical efficiency of our approach, we compare our algorithm with the best-known heuristics. Notably, these state-of-the-art heuristics hold a priori knowledge of specific problems to tailor their strategies, thereby exploiting the inherent characteristics of these problems to improve efficiency. In contrast, our generalized approach increases the challenge by omitting the prior knowledge and problem features. However, this approach offers the advantage of adaptability, eliminating the need to redesign in the context of new problems.

For the experiments, we employ GUROBI version 9.5.1 as a branch-and-cut (B\&C) solver, operating in a black-box way without customized algorithmic configurations. We set the sample size $N$ to 6000 and the parameters $\pi_{card} = 3$, $\pi_{limit}=5$, and $R_{\text{limit}} = 0.5$ in Algorithm \ref{ag: train}. The parameter $\theta$ in Algorithm \ref{ag: train} is set to 0.2 in Section \ref{subsec: expcase1}, and in Section \ref{subsec: expcase2}, it is 0.2 if $|J|<40$ and 0.5 otherwise. The parameter $\theta$ in Algorithm \ref{ag: train} is set to 0.2 in \ref{subsec: expcase1} and The setting of the regulation parameter in the regression follows \cite{zou2007degrees}. Experiments are implemented in the Python language. We use a workstation of Intel Xeon E5-2699 v4 2.20 gigahertz CPU for sampling and a personal computer (Intel i7 2.60 gigahertz CPU, 16-gigabyte memory) with Windows 10.0. for comparison tests. Additionally, to facilitate a clear and structured presentation of our comparative analysis, we introduce specific notations in Table \ref{tab: notations}. These notations are designed to display the relative performance and efficiency data between the tested models.
\begin{table}[]
\TABLE
{{Parameters for instance generation}\label{tab: notations}}
{\resizebox{\columnwidth}{!}{%
\begin{tabularx}{500pt}{lp{420pt}}
\hline 
Notation & Description  \\
\hline SIM &   The set indicator model\\
Exact & The exact model for comparison  \\
$\Delta$ & Instance group identifier \\
Nodes& Number of explored nodes within the B\&C procedure \\ 
Rows& Number of constraints of a MIP before presolving\\
Cols& Number of variables of a MIP before presolving\\
Obj&  Objective value (computed by equation (\ref{eq: original}) for SIM)\\
$\text{Time}_{\text{MIP}}$ & CPU computing time in seconds of solving MIP models of SIM or Exact\\
$\text{Time}_{\text{Total}}$& Total CPU computing time in seconds, including the average time to conduct Algorithm \ref{ag: Monte-Carlo} once and the time for Algorithm \ref{ag: train}.\\
Gap& Integrality gap\\
Alg1 / Alg2& Objective ratio of algorithms denoted by Alg1 and Alg2.\\
Min(sample)& Minimal objective value of samples\\
\hline
\end{tabularx}
}
}
{}
\end{table}
\subsection{Case study I: the joint order batching and picker routing problem (JOBPRP)}
\label{subsec: expcase1}

We adopt the problem setting of \cite{valle2017optimally} to JOBPRP with a warehouse configuration that comprises multiple identical vertical storage racks, separated by navigable aisles. Each rack is divided into several blocks, with each block storing a specific type of item. Orders are composed of one or more items, and items belonging to the same order must be grouped into a single batch. The process of handling a batch involves a picker departing from a depot, traversing through the blocks for picking, and subsequently returning to the depot. The distance between any two blocks is calculated using the Manhattan distance metric.

The JOBPRP has attracted significant academic interest in the field of warehouse management studies, as highlighted in the comprehensive review conducted by \cite{van2018designing}. The JOBPRP is initially proposed by \cite{won2005joint}. For exact B\&C algorithms, the work by \cite{valle2016modelling} presented three basic formulations. Two of them are based on network flows, with the remaining one having exponentially increasing constraints. Furthermore, real-world data sets were included in the study of this problem \citep{valle2016modelling}. In a subsequent advancement, the work by \cite{valle2017optimally} enhanced these formulations with valid inequalities. The study by \cite{zhang2023improved} further made significant improvements by reformulating the connectivity constraints within the JOBPRP. Additionally, several heuristics and meta-heuristics have also been proposed for this problem in the past decade \cite{cheng2015using, lin2016joint, li2017joint, valle2020order}.

%\cite{cheng2015using, lin2016joint, li2017joint, briant2020efficient}

In this study, we compare with the latest study of the exact approach by \cite{zhang2023improved}. Compared with \cite{valle2017optimally}, the new formulation by \cite{zhang2023improved} increases the number of optimally solved instances within 2400 seconds by 50\%. In scenarios involving over ten orders where both algorithms achieve optimal solutions, the new formulation ensures a minimum of 54\% reduction in solution time. We adopt the same benchmark instances as \cite{valle2017optimally, zhang2023improved}. 

\subsubsection{Experiment results}
Table \ref{tab: JOBPRP} presents detailed information and solutions on 39 instances from the SIM and Exact model. On average, the distance in objective values between SIM and the Exact is within 1\%.  The SIM performs equal or better in one-third of all the instances. On average, SIM reduces the solution time by 95\%,  demonstrating a significant total time advantage. These results strongly indicate the precision and efficiency of the SIM, as answers to questions \textbf{Q1} and \textbf{Q2}.

\begin{figure}
	\centering
	\begin{minipage}{.35\textwidth}
		\FIGURE
		{\includegraphics[scale=0.35]{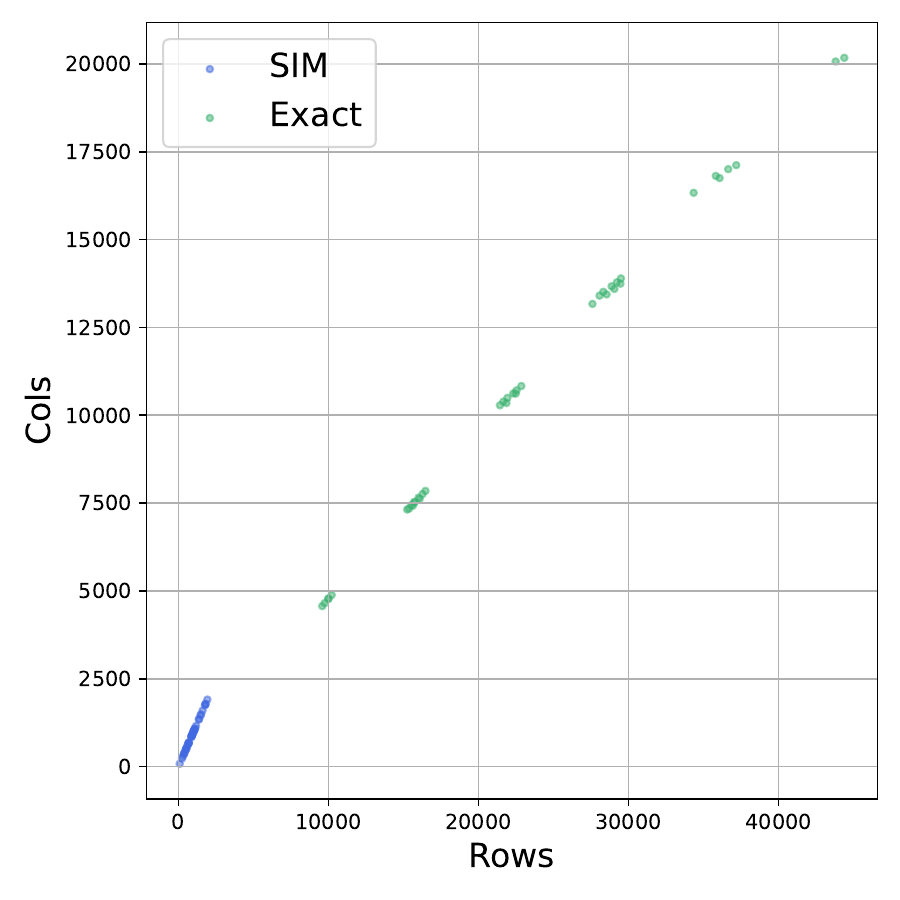}}
      %\captionsetup{width=0.9\textwidth}
		{{Scatter of model sizes}\label{fig: OBPscatter}}
        {}
	\end{minipage}%
	\begin{minipage}{.65\textwidth}
		\FIGURE
        {\includegraphics[scale=0.35]{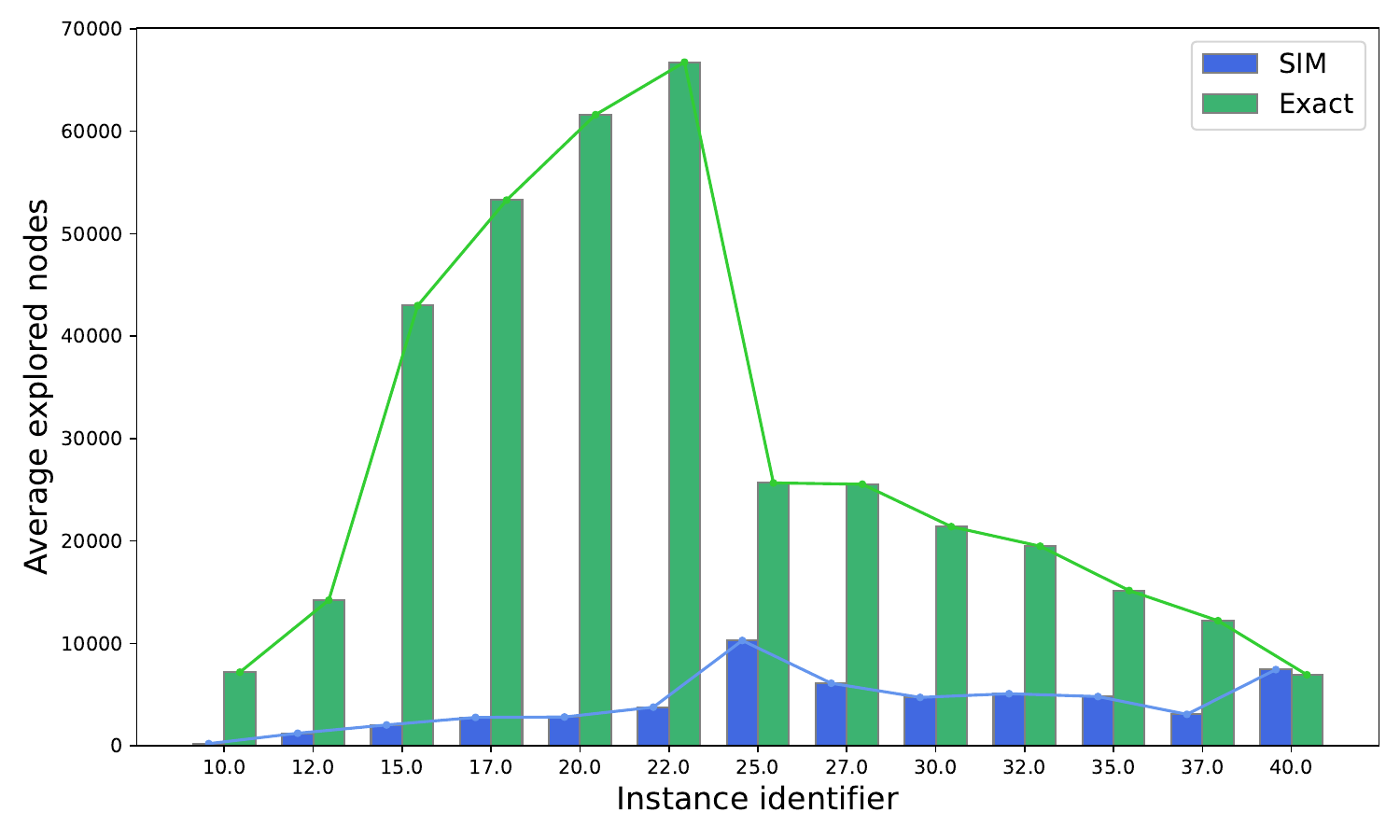}}
            %\captionsetup{width=0.9\textwidth}
		{{Explored nodes during the B\&C process}\label{fig: OBPnodes}}
        {}
	\end{minipage}
\end{figure}

We visualize part of the data for further analysis. Figure \ref{fig: OBPscatter} depicts a scatter plot of the number of constraints versus the number of variables for both the SIM and the Exact across different instances, with the number of constraints on the x-axis and the number of variables on the y-axis. It is observable that the distribution of the SIM is relatively concentrated within the range of 80-2000. In contrast, the distribution of data for the Exact shows an incremental trend with the increase in the scale of instances. This demonstrates that our framework is capable of stably generating MIP models of far less scale. 

Figure \ref{fig: OBPnodes} presents a bar graph illustrating the average explored nodes during the B\&C process for different amounts of orders. For the Exact, the node number increases rapidly before 22 orders, caused by the rapid growth of model size. Beyond 22 orders, due to the time limit and model complexity, the B\&C process stops early, leaving a sharp decline followed by a gradual decrease in the figure. Conversely, for the SIM, the nodes remain at a lower level with a slow growth trend. This aligns with the SIM‘s limited number of variables in Figure \ref{fig: OBPscatter}.

\begin{figure}
	\centering
	\begin{minipage}{.35\textwidth}
		\FIGURE
       %\captionsetup{width=0.9\textwidth}
      {\includegraphics[scale=0.24]{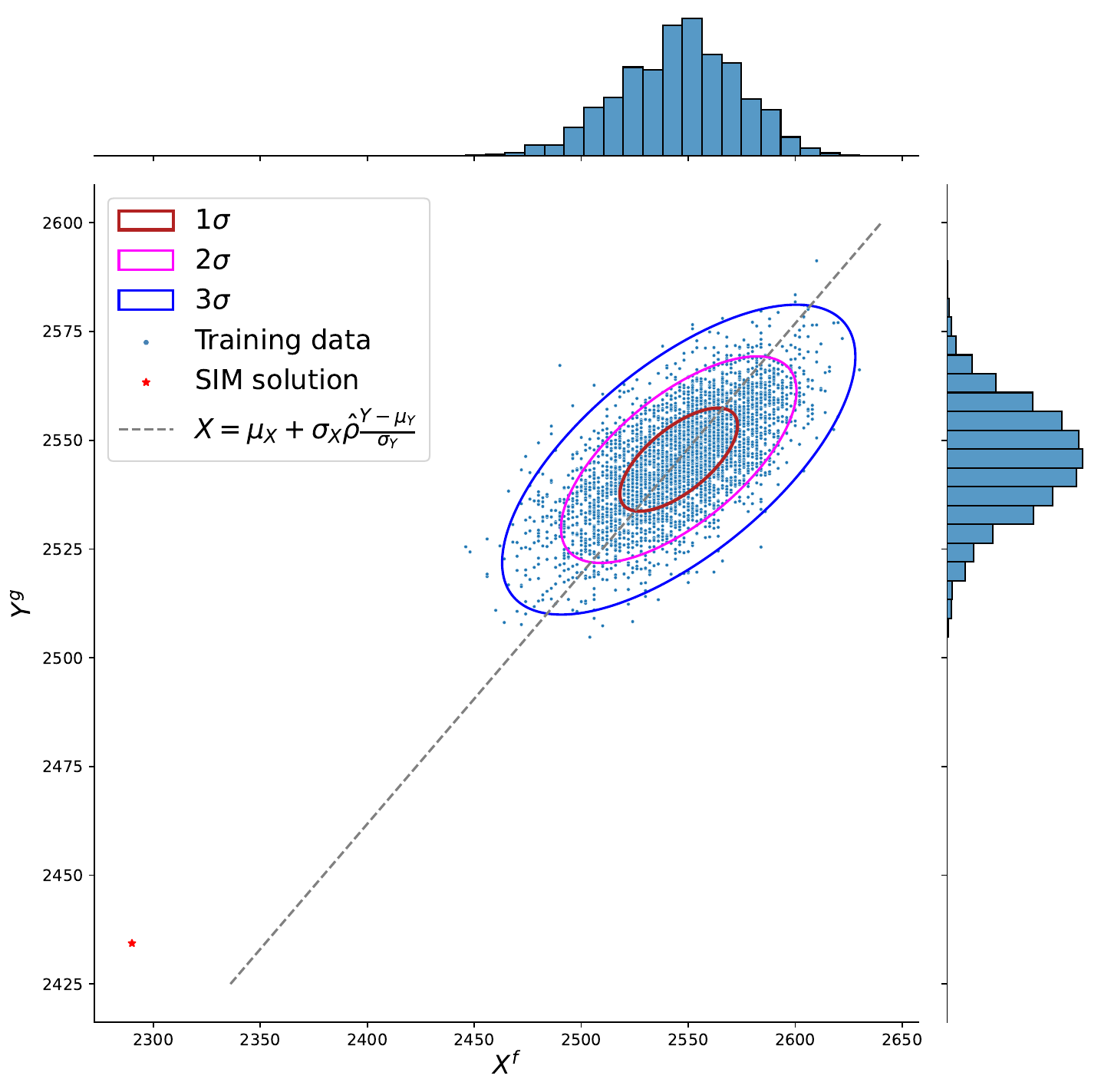}}
      %\captionsetu{width=.9\linewidth}
		{{Parity plot of training and prediction of the SIM and Exact models}
		\label{fig: OBPparity}}
        {}
	\end{minipage}%
	\begin{minipage}{.65\textwidth}
		\FIGURE
        {\includegraphics[scale=0.33]{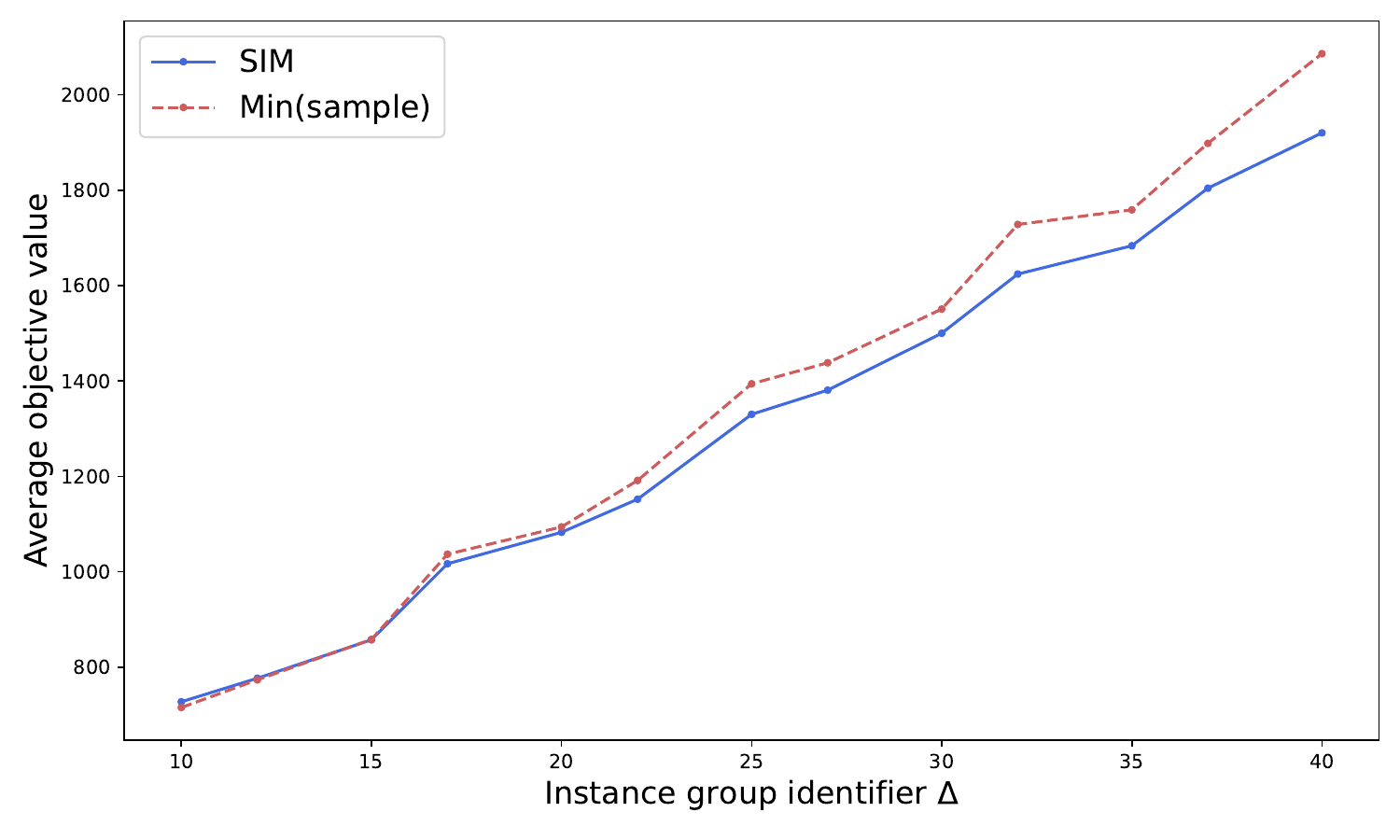}}
            %\captionsetup{width=0.9\textwidth}
        %\captionsetup{width=.9\linewidth}
		{{The comparison of average group objective values of the SIM the Min(sample)}
          \label{fig: OBPsample}}
        {}
	\end{minipage}
\end{figure}

In Figure \ref{fig: OBPparity}, we visualize the collected data for the instances of Day20 with 40 orders, as a typical example. The scatter plot presents the Exact objective values on the x-axis and the learned approximate function g values on the y-axis. This plot exhibits a bivariate normal distribution. The formula (\ref{eq: expectation}) from Theorem \ref{thm: bivariate} is represented as a straight line through statistical estimations. It is observed that the solution by the SIM is situated close to its expected value and is significantly lower than all sampled values. In Figure \ref{fig: OBPsample}, we further illustrate the average SIM objective values for different amounts of orders alongside the minimum sampled objective values. The objective values from the SIM are increasingly lower than that of the sampled values.  This evidence indicates that our framework outperforms mere Monte Carlo sampling, being able to learn optimization models to discover new solutions.

The remarkable improvements demonstrated by the SIM can be attributed to our modifications of the MIP formulations, specifically, the replacement of routing constraints with big-M constraints. Thus the question \textbf{Q3} is answered. Firstly, considering that the experiment setting utilizes the same solving procedure, the only varying factor is the formulation structure. Moreover, the empirical evidence shown in Figure \ref{fig: OBPscatter} and \ref{fig: OBPnodes} indicates that our framework consistently generates MIP models of smaller size, alongside reductions in the B\&B tree sizes. These outcomes corroborate the efficiency shown in the experiment and further prove our motivation expressed in the introduction. Finally, the significant difference in Figure \ref{fig: OBPsample} between the SIM and samples demonstrates that our framework does have the ability to find optimization direction and find new solutions. 

Additionally, we compare our approach against the best-known heuristic for JOBPRP. We imposed a 20-second time restriction on each surrogate model within Algorithm \ref{ag: train}, recognizing that the solution at this limit could represent the optimal outcome. \cite{valle2020order} designates this best-known heuristic as PIO, reporting its superior performance over traditional heuristics in effectiveness and efficiency. We utilized the same problem scales as \cite{valle2020order}. The comparative results are presented in Table \ref{tab: JOBPRPH}. The results indicate that our algorithm consistently performs better on small and medium scales. On average, our approach demonstrates better outcomes and reduced computational time.

\begin{table}[htbp]
  \TABLE
  {{Results for the JOBPRP instances}  \label{tab: JOBPRP}}
  {\resizebox{\columnwidth}{!}{%
    \begin{tabular}{rrr|rrrrrrrrrrrrrr}
    \toprule
          \multicolumn{3}{c}{Instance}  & \multicolumn{5}{c}{SIM} & \multicolumn{6}{c}{Exact}                       &       &  \\
          \cmidrule(lr){1-3}
          \cmidrule(lr){4-9}
          \cmidrule(lr){10-15}
          $\Delta$ & Order & \multicolumn{1}{c}{Trolley} & Obj & Nodes&Rows & Cols & $\text{Time}_{\text{MIP}}$ & $\text{Time}_{\text{Total}}$ &Obj& Nodes & Rows & Cols &Gap(\%)&$\text{Time}_{\text{MIP}}$ & Exact/SIM & Min(sample) \\
          \hline
          & 10    & 2     & 582   & 145   & 290   & 288   & 0.2   & 18    & 578   & 3040  & 9582  & 4568  & 0.00  & 10    & 0.99  & 578 \\
          & 12    & 2     & 630   & 332   & 596   & 594   & 0.8   & 14    & 616   & 4696  & 9742  & 4652  & 0.00  & 73    & 0.98  & 616 \\
          & 15    & 2     & 652   & 1     & 544   & 542   & 0.4   & 23    & 650   & 3953  & 9992  & 4786  & 0.00  & 52    & 1.00  & 650 \\
          & 17    & 3     & 832   & 4316  & 682   & 678   & 5.9   & 61    & 802   & 15481 & 15244 & 7317  & 0.00  & 271   & 0.96  & 864 \\
          & 20    & 3     & 880   & 221   & 364   & 360   & 0.3   & 24    & 864   & 50341 & 15529 & 7446  & 0.00  & 917   & 0.98  & 886 \\
          & 22    & 3     & 892   & 2497  & 1006  & 1002  & 9.5   & 130   & 892   & 107936 & 15724 & 7536  & 0.00  & 1923  & 1.00  & 946 \\
    \multicolumn{1}{l}{Day5} & 25    & 4     & 1120  & 7266  & 999   & 992   & 13.4  & 83    & 1108  & 25543 & 21427 & 10284 & 25.72  & 2400  & 0.99  & 1190 \\
          & 27    & 4     & 1180  & 4367  & 891   & 884   & 15.1  & 90    & 1154  & 25575 & 21651 & 10388 & 25.74  & 2400  & 0.98  & 1248 \\
          & 30    & 4     & 1242  & 4539  & 695   & 688   & 6.0   & 37    & 1228  & 25762 & 21927 & 10496 & 29.64  & 2400  & 0.99  & 1292 \\
          & 32    & 5     & 1358  & 3575  & 1046  & 1035  & 5.8   & 26    & 1346  & 24390 & 27591 & 13170 & 38.48  & 2400  & 0.99  & 1504 \\
          & 35    & 5     & 1464  & 3929  & 1171  & 1160  & 14.5  & 43    & 1448  & 13681 & 28061 & 13405 & 42.33  & 2400  & 0.99  & 1548 \\
          & 37    & 5     & 1538  & 4750  & 846   & 835   & 6.4   & 39    & 1506  & 15362 & 28316 & 13515 & 41.50  & 2400  & 0.98  & 1624 \\
          & 40    & 6     & 1712  & 15454 & 1924  & 1908  & 67.9  & 183   & 1642  & 8612  & 34336 & 16332 & 49.90  & 2400  & 0.96  & 1848 \\
          \hline
           & 10    & 2     & 656   & 1     & 82    & 80    & 0.1   & 9     & 656   & 7598  & 9988  & 4768  & 0.00  & 34    & 1.00  & 656 \\
          & 12    & 2     & 702   & 1     & 490   & 488   & 0.5   & 13    & 702   & 2501  & 10212 & 4884  & 0.00  & 40    & 1.00  & 702 \\
          & 15    & 3     & 894   & 4326  & 862   & 858   & 5.6   & 24    & 874   & 39134 & 15727 & 7515  & 0.00  & 315   & 0.98  & 888 \\
          & 17    & 3     & 968   & 1072  & 487   & 483   & 0.9   & 16    & 960   & 76979 & 16000 & 7653  & 0.00  & 703   & 0.99  & 968 \\
          & 20    & 3     & 1010  & 372   & 373   & 369   & 0.2   & 14    & 984   & 104755 & 16267 & 7758  & 0.00  & 1523  & 0.97  & 1026 \\
          & 22    & 3     & 1022  & 762   & 403   & 399   & 0.4   & 16    & 1000  & 66779 & 16462 & 7848  & 0.00  & 1109  & 0.98  & 1032 \\
    \multicolumn{1}{l}{Day10} & 25    & 4     & 1236  & 18756 & 1603  & 1596  & 91.9  & 225   & 1196  & 25603 & 22311 & 10620 & 26.84  & 2400  & 0.97  & 1294 \\
          & 27    & 4     & 1274  & 6731  & 963   & 956   & 6.4   & 28    & 1274  & 25219 & 22531 & 10708 & 26.50  & 2400  & 1.00  & 1350 \\
          & 30    & 4     & 1304  & 1184  & 903   & 896   & 4.2   & 45    & 1286  & 25145 & 22851 & 10832 & 31.25  & 2400  & 0.99  & 1360 \\
          & 32    & 5     & 1548  & 4014  & 1501  & 1490  & 32.3  & 125   & 1492  & 20177 & 28866 & 13670 & 37.33  & 2400  & 0.96  & 1628 \\
          & 35    & 5     & 1584  & 4011  & 1016  & 1005  & 8.1   & 35    & 1558  & 19281 & 29216 & 13785 & 40.44  & 2400  & 0.98  & 1642 \\
          & 37    & 5     & 1648  & 361   & 926   & 915   & 1.6   & 24    & 1590  & 15872 & 29491 & 13895 & 39.18  & 2400  & 0.96  & 1706 \\
          & 40    & 6     & 1784  & 3713  & 1804  & 1788  & 45.9  & 133   & 1800  & 6874  & 35812 & 16812 & 55.28  & 2400  & 1.01  & 1972 \\
          \hline
          & 10    & 3     & 944   & 491   & 232   & 228   & 0.3   & 11    & 912   & 10920 & 15367 & 7344  & 0.00  & 112   & 0.97  & 912 \\
          & 12    & 3     & 998   & 3350  & 646   & 642   & 2.1   & 17    & 998   & 35459 & 15613 & 7422  & 0.00  & 254   & 1.00  & 1002 \\
          & 15    & 3     & 1026  & 1787  & 505   & 501   & 1.3   & 17    & 1022  & 85837 & 16078 & 7623  & 0.00  & 532   & 1.00  & 1036 \\
          & 17    & 4     & 1250  & 2929  & 675   & 668   & 4.9   & 23    & 1250  & 67364 & 21859 & 10348 & 5.84  & 2400  & 1.00  & 1278 \\
          & 20    & 4     & 1358  & 7842  & 871   & 864   & 26.1  & 90    & 1348  & 29741 & 22479 & 10616 & 15.87  & 2400  & 0.99  & 1370 \\
          & 22    & 5     & 1542  & 8062  & 1091  & 1080  & 41.2  & 127   & 1570  & 25508 & 28531 & 13440 & 25.74  & 2400  & 1.02  & 1596 \\
    \multicolumn{1}{l}{Day20} & 25    & 5     & 1634  & 4828  & 1091  & 1080  & 39.5  & 236   & 1670  & 25842 & 29051 & 13595 & 34.43  & 2400  & 1.02  & 1698 \\
          & 27    & 5     & 1688  & 7221  & 1066  & 1055  & 40.4  & 159   & 1700  & 25824 & 29456 & 13745 & 32.82  & 2400  & 1.01  & 1716 \\
          & 30    & 6     & 1954  & 8451  & 1354  & 1338  & 22.7  & 80    & 1940  & 13284 & 36064 & 16752 & 47.42  & 2400  & 0.99  & 2000 \\
          & 32    & 6     & 1966  & 7668  & 1372  & 1356  & 19.7  & 114   & 1900  & 13888 & 36640 & 17004 & 50.85  & 2400  & 0.97  & 2052 \\
          & 35    & 6     & 2002  & 6482  & 1474  & 1458  & 32.8  & 95    & 2030  & 12516 & 37174 & 17118 & 46.15  & 2400  & 1.01  & 2086 \\
          & 37    & 7     & 2226  & 4110  & 1765  & 1743  & 34.0  & 182   & 2246  & 5403  & 43800 & 20069 & 57.20  & 2400  & 1.01  & 2364 \\
          & 40    & 7     & 2264  & 3150  & 1793  & 1771  & 47.1  & 315   & 2330  & 5306  & 44374 & 20174 & 59.80  & 2400  & 1.03  & 2438 \\
          \hline
         \multicolumn{3}{c}{Average} &  & 4181.2  & 933.4  &924.9       & 16.8  & 75.4  &       & 28645.7  & 23367.7  & 11022.9  &       & 1678.6  & 0.99  &  \\
          \bottomrule
    \end{tabular}%
    }
    }
    {}
\end{table}%

\begin{table}[htbp]
  \TABLE
  {{The SIM in JOBPRP heuristic comparison}  \label{tab: JOBPRPH}}%
  %\resizebox{\columnwidth}{!}{%
  {\begin{tabular}{rrrrrrrr}
    \toprule
    \multicolumn{3}{c}{Instance} & \multicolumn{2}{c}{SIM(20s)} & \multicolumn{2}{c}{PIO} &\\
          \cmidrule(lr){1-3}
          \cmidrule(lr){4-5}
          \cmidrule(lr){6-7}
        $\Delta$  & \multicolumn{1}{l}{Order} & \multicolumn{1}{l}{Trolley} & \multicolumn{1}{l}{Obj} & \multicolumn{1}{l}{$\text{Time}_{\text{Total}}$} & \multicolumn{1}{l}{Obj} & \multicolumn{1}{l}{$\text{Time}_{\text{Total}}$} &  PIO/SIM\\
    \toprule
    \multicolumn{1}{l}{Day5} & 25    &       & \textbf{1120}  & 67    & 1137.8 & 1     & 1.02  \\
          & 30    &       & \textbf{1198}  & 167   & 1267.9 & 3     & 1.06  \\
          & 50    &       & 2032 & 154  & \textbf{2010.5} & 23    & 0.99  \\
          & 75    &       & 2974  & 297   & \textbf{2884.8} & 648   & 0.97  \\
          \hline
    \multicolumn{1}{l}{Day10} & 25    &       & \textbf{1236}  & 34    & 1251.3 & 1     & 1.01  \\
          & 30    &       & \textbf{1300}  & 26    & 1364.9 & 3     & 1.05  \\
          & 50    &       & \textbf{2174}  & 175   & 2184.4 & 29    & 1.00  \\
          & 75    &       & 3140  & 776   & \textbf{3125.8} & 1903  & 1.00  \\
          \hline
    \multicolumn{1}{l}{Day20} & 25    &       & \textbf{1630}  & 102   & 1708.8 & 2     & 1.05  \\
          & 30    &       & \textbf{1908}  & 149   & 1954.8 & 3     & 1.02  \\
          & 50    &       & \textbf{2688}  & 142   & 2754.6 & 54    & 1.02  \\
          & 75    &       & 4006  & 1120  & \textbf{3849.0}  & 659   & 0.96  \\
          \hline
          \multicolumn{3}{c}{Average}  &     & 267.4  &       & 277.4  & 1.01  \\
    \bottomrule
    \end{tabular}%
    }
    {}
\end{table}%

\subsection{Case Study II: the soft-clustered vehicle routing problem (SoftCluVRP)}
\label{subsec: expcase2}
The SoftCluVRP represents a variant of the classical CVRP. In SoftCluVRP, customers are partitioned into distinct clusters, with the restriction that all customers belonging to the same cluster must be served by the same vehicle. This requirement introduces an additional layer of complexity compared to the classical CVRP. Moreover, unlike the hard-clustered variant of this problem, where a vehicle must completely serve one cluster before proceeding to the next, SoftCluVRP allows for a more flexible approach. In this variant, a vehicle can interrupt service to a cluster to serve customers in another cluster, and then return to the previous cluster if necessary. This flexibility in service sequence allows more efficient routing strategies. However, this flexibility also increases the complexity involved in routing decisions. The problem aims to find an optimal set of routes for a fleet of vehicles that minimizes the total travel cost while ensuring that each cluster's customers are exclusively served by a single vehicle, and adhering to classical CVRP constraints. Moreover, to enhance the diversity of $\Bar{H}$ in Algorithm \ref{ag: train} for large instances, we apply a random permutation each time it is newly assigned. As our algorithm becomes non-deterministic in this scenario, we perform 5 experimental rounds per instance, recording the best and average objective values as $\text{Obj}_{\text{best}}$ and $\text{Obj}_{\text{avg}}$, respectively. Additional model data for the SIM is also displayed along with their average values.

The clustered vehicle routing problem was initially introduced by \cite{sevaux2008hamiltonian} as a practical problem. Subsequently, \cite{defryn2017fast} involved soft constraints in customer cluster assignments, implying that such modifications could result in significant cost reduction. Exact methods have been proposed by \cite{hintsch2020exact}. Heuristics have also been proposed by \cite{hintsch2021large}. We refer to the review in the study by \cite{hintsch2021large} for a detailed investigation. In this experiment, we adopt the latest B\&C algorithm by \cite{hessler2021branch}, which has been examined to be significantly effective. We use the well-known \texttt{Golden} benchmark instances \citep{hintsch2020exact, hintsch2021large, hessler2021branch}. This benchmark encompasses a collection of 220 large-scale instances, with 201 to 484 customers and 14 to 97 clusters. The benchmark is derived from 20 original instances, distinguished by varying clusters, and grouped based on the identifiers of the original instances. Limited by computing resources to conduct Monte-Carlo on large instances, we conduct our comparison on 110 instances with customers no more than 400.  

\begin{table}[htbp]
  \TABLE
  {{Detailed results for the SoftCluVRP instances (Part1) }  \label{tab: SoftCluE1}}
   {\resizebox{\columnwidth}{!}{%
    \begin{tabular}{lllrrrrrrrrrrrrrrrr}
        \toprule
        \multicolumn{1}{c}{\multirow{2}[0]{*}{$\Delta$}}  &\multicolumn{1}{c}{\multirow{2}[0]{*}{Cluster}} &\multicolumn{1}{c}{\multirow{2}[0]{*}{Customer}}  & \multicolumn{8}{c}{SIM} & \multicolumn{6}{c}{Exact}& &  \\
        \cmidrule(lr){4-11}
        \cmidrule(lr){12-17}
        \cmidrule(lr){18-19}
        & & &$\text{Obj}_{\text{best}}$ &$\text{Obj}_{\text{avg}}$ &\multicolumn{1}{l}{Nodes} & \multicolumn{1}{l}{Rows} & \multicolumn{1}{l}{Cols} & \multicolumn{1}{l}{Gap($\%$)} & $\text{Time}_{\text{MIP}}$ & $\text{Time}_{\text{Total}}$& Obj&Nodes &Rows &Cols & Gap($\%$) &$\text{Time}_{\text{MIP}}$& \multicolumn{1}{l}{Min(sample)}& $\text{Exact/SIM}_{\text{best}}$ \\
     \hline
    1     & 17    & 241   & 4640  & 4640.0  & 1234  & 2345  & 2388  & 0.0   & 2     & 6     & \multicolumn{1}{r}{4640} & 465888 & 28658 & 29073 & \multicolumn{1}{r}{1.1 } & 3600  & 4868  & \multicolumn{1}{r}{1.00 } \\
    1     & 18    & 241   & 4652  & 4664.6  & 3780  & 3376  & 3422  & 0.0   & 5     & 9     & \multicolumn{1}{r}{4968} & 220456 & 29676 & 29091 & \multicolumn{1}{r}{10.5 } & 3600  & 4884  & \multicolumn{1}{r}{1.07 } \\
    1     & 19    & 241   & 4667  & 4685.0  & 3532  & 3406  & 3455  & 0.0   & 2     & 10    & \multicolumn{1}{r}{4791} & 241965 & 30207 & 29110 & \multicolumn{1}{r}{6.1 } & 3600  & 4982  & \multicolumn{1}{r}{1.03 } \\
    1     & 21    & 241   & 4669  & 4681.0  & 6281  & 3402  & 3457  & 0.0   & 9     & 13    & \multicolumn{1}{r}{5586} & 198212 & 31380 & 29151 & \multicolumn{1}{r}{20.5 } & 3600  & 5073  & \multicolumn{1}{r}{1.20 } \\
    1     & 22    & 241   & 4677  & 4694.0  & 3490  & 3433  & 3491  & 0.0   & 7     & 25    & \multicolumn{1}{r}{5027} & 213993 & 32087 & 29173 & \multicolumn{1}{r}{11.7 } & 3600  & 5164  & \multicolumn{1}{r}{1.07 } \\
    1     & 25    & 241   & 4660  & 4664.6  & 3803  & 3368  & 3435  & 0.0   & 7     & 30    & \multicolumn{1}{r}{4932} & 165596 & 34467 & 29245 & \multicolumn{1}{r}{98.3 } & 3600  & 5349  & \multicolumn{1}{r}{1.06 } \\
    1     & 27    & 241   & 4658  & 4663.4  & 3044  & 3495  & 3568  & 0.0   & 8     & 23    & \multicolumn{1}{r}{5572} & 214888 & 36612 & 29298 & \multicolumn{1}{r}{20.4 } & 3600  & 5415  & \multicolumn{1}{r}{1.20 } \\
    1     & 31    & 241   & 4703  & 4710.2  & 10098  & 3238  & 3323  & 0.0   & 30    & 108   & \multicolumn{1}{r}{4941} & 65436 & 41727 & 29416 & \multicolumn{1}{r}{10.7 } & 3600  & 5532  & \multicolumn{1}{r}{1.05 } \\
    1     & 35    & 241   & 4712  & 4731.0  & 12453  & 3542  & 3639  & 0.0   & 36    & 99    & \multicolumn{1}{r}{6102} & 49055 & 48235 & 29550 & \multicolumn{1}{r}{28.4 } & 3600  & 5624  & \multicolumn{1}{r}{1.29 } \\
    1     & 41    & 241   & 4687  & 4757.0  & 8339  & 9782  & 9897  & 0.0   & 28    & 127   & -     & 6458  & 60927 & 29781 & -     & 3600  & 5595  & - \\
    1     & 49    & 241   & 4803  & 4902.2  & 4087  & 9669  & 9808  & 0.0   & 41    & 157   & \multicolumn{1}{r}{4785} & 214   & 84712 & 30145 & \multicolumn{1}{r}{8.2 } & 3600  & 5828  & \multicolumn{1}{r}{1.00 } \\
    \\
    2     & 22    & 321   & 7442  & 7482.2  & 8014  & 3405  & 3463  & 0.0   & 12    & 28    & \multicolumn{1}{r}{8684} & 65931 & 52875 & 51613 & \multicolumn{1}{r}{19.1 } & 3600  & 7998  & \multicolumn{1}{r}{1.17 } \\
    2     & 23    & 321   & 7392  & 7486.2  & 18081  & 3316  & 3377  & 0.0   & 16    & 25    & -     & 112507 & 53633 & 51636 & -     & 3600  & 7932  & - \\
    2     & 25    & 321   & 7392  & 7429.2  & 5276  & 3388  & 3455  & 0.0   & 13    & 44    & -     & 49164 & 55350 & 51685 & -     & 3600  & 8038  & - \\
    2     & 27    & 321   & 7399  & 7448.0  & 8545  & 3286  & 3359  & 0.0   & 22    & 32    & \multicolumn{1}{r}{7949} & 57483 & 57769 & 51738 & \multicolumn{1}{r}{11.6 } & 3600  & 8132  & \multicolumn{1}{r}{1.07 } \\
    2     & 30    & 321   & 7440  & 7446.0  & 21396  & 3430  & 3512  & 0.0   & 41    & 81    & -     & 26652 & 61805 & 51825 & -     & 3600  & 8246  & - \\
    2     & 33    & 321   & 7359  & 7460.0  & 3235  & 3421  & 3512  & 0.0   & 25    & 100   & -     & 6643  & 66362 & 51921 & -     & 3600  & 8225  & - \\
    2     & 36    & 321   & 7359  & 7445.0  & 2406  & 3824  & 3924  & 0.0   & 7     & 100   & -     & 6796  & 71653 & 52026 & -     & 3600  & 8371  & - \\
    2     & 41    & 321   & 7408  & 7428.0  & 1768  & 9329  & 9444  & 0.0   & 49    & 130   & -     & 6760  & 82597 & 52221 & -     & 3600  & 8530  & - \\
    2     & 46    & 321   & 7449  & 7500.0  & 24705  & 9790  & 9920  & 0.0   & 85    & 263   & -     & 645   & 96640 & 52441 & -     & 3600  & 8658  & - \\
    2     & 54    & 321   & 7591  & 7720.0  & 25428  & 9838  & 9992  & 0.0   & 124   & 259   & -     & 177   & 125988 & 52845 & -     & 3600  & 9103  & - \\
    2     & 65    & 321   & 7745  & 7816.0  & 6296  & 9811  & 9998  & 0.0   & 30    & 324   & -     & 6     & 183430 & 53505 & -     & 3600  & 9470  & - \\
    \\
    5     & 14    & 201   & 6970  & 6970.0  & 1069  & 3289  & 3323  & 0.0   & 2     & 19    & \multicolumn{1}{r}{7084} & 161712 & 18925 & 20205 & \multicolumn{1}{r}{3.1 } & 3600  & 7231  & \multicolumn{1}{r}{1.02 } \\
    5     & 15    & 201   & 6933  & 6933.0  & 684   & 2402  & 2426  & 0.0   & 2     & 5     & \multicolumn{1}{r}{6742} & 58673 & 19295 & 20220 & \multicolumn{1}{r}{0.0 } & 2118  & 6989  & \multicolumn{1}{r}{0.97 } \\
    5     & 16    & 201   & 6895  & 6920.2  & 1426  & 2371  & 2397  & 0.0   & 2     & 10    & \multicolumn{1}{r}{6742} & 40073 & 19704 & 20236 & \multicolumn{1}{r}{0.0 } & 1777  & 7099  & \multicolumn{1}{r}{0.98 } \\
    5     & 17    & 201   & 6862  & 6862.0  & 908   & 2466  & 2494  & 0.0   & 2     & 5     & \multicolumn{1}{r}{7118} & 101846 & 20180 & 20253 & \multicolumn{1}{r}{5.0 } & 3600  & 7377  & \multicolumn{1}{r}{1.04 } \\
    5     & 19    & 201   & 7030  & 7030.0  & 3758  & 3321  & 3370  & 0.0   & 8     & 28    & \multicolumn{1}{r}{6874} & 14062 & 21611 & 20290 & \multicolumn{1}{r}{0.0 } & 993   & 8019  & \multicolumn{1}{r}{0.98 } \\
    5     & 21    & 201   & 6881  & 6885.0  & 6805  & 3382  & 3437  & 0.0   & 10    & 20    & \multicolumn{1}{r}{6887} & 80144 & 22937 & 20331 & \multicolumn{1}{r}{2.1 } & 3600  & 7796  & \multicolumn{1}{r}{1.00 } \\
    5     & 23    & 201   & 6887  & 6943.8  & 5235  & 3431  & 3492  & 0.0   & 14    & 23    & \multicolumn{1}{r}{6750} & 17155 & 24370 & 20376 & \multicolumn{1}{r}{0.0 } & 548   & 7823  & \multicolumn{1}{r}{0.98 } \\
    5     & 26    & 201   & 6880  & 6882.6  & 8886  & 3746  & 3816  & 0.0   & 7     & 18    & \multicolumn{1}{r}{6704} & 25846 & 26975 & 20451 & \multicolumn{1}{r}{0.0 } & 585   & 8420  & \multicolumn{1}{r}{0.97 } \\
    5     & 29    & 201   & 6879  & 6910.0  & 3861  & 3737  & 3816  & 0.0   & 17    & 25    & \multicolumn{1}{r}{6704} & 160440 & 30239 & 20535 & \multicolumn{1}{r}{0.0 } & 3161  & 8243  & \multicolumn{1}{r}{0.97 } \\
    5     & 34    & 201   & 6833  & 6900.0  & 6035  & 3826  & 3920  & 0.0   & 17    & 35    & \multicolumn{1}{r}{6684} & 58033 & 37572 & 20695 & \multicolumn{1}{r}{1.8 } & 3600  & 8627  & \multicolumn{1}{r}{0.98 } \\
    5     & 41    & 201   & 6752  & 6752.0  & 441   & 9935  & 10050  & 0.0   & 2     & 74    & \multicolumn{1}{r}{6727} & 17748 & 52304 & 20961 & \multicolumn{1}{r}{4.4 } & 3600  & 8781  & \multicolumn{1}{r}{1.00 } \\
    \\
    6     & 19    & 281   & 8115  & 8122.2  & 973   & 2263  & 2295  & 0.0   & 3     & 14    & \multicolumn{1}{r}{8115} & 212140 & 38268 & 39530 & \multicolumn{1}{r}{0.0 } & 3180  & 8593  & \multicolumn{1}{r}{1.00 } \\
    6     & 21    & 281   & 8119  & 8134.4  & 1786  & 2413  & 2449  & 0.0   & 3     & 7     & \multicolumn{1}{r}{8799} & 194659 & 39837 & 39571 & \multicolumn{1}{r}{9.7 } & 3600  & 8672  & \multicolumn{1}{r}{1.08 } \\
    6     & 22    & 281   & 8107  & 8107.0  & 4627  & 2384  & 2422  & 0.0   & 6     & 15    & \multicolumn{1}{r}{8841} & 172221 & 41062 & 39593 & \multicolumn{1}{r}{10.3 } & 3600  & 8730  & \multicolumn{1}{r}{1.09 } \\
    6     & 24    & 281   & 8316  & 8351.0  & 1197  & 3385  & 3449  & 0.0   & 3     & 8     & \multicolumn{1}{r}{9035} & 157285 & 43035 & 39640 & \multicolumn{1}{r}{10.6 } & 3600  & 9401  & \multicolumn{1}{r}{1.09 } \\
    6     & 26    & 281   & 8285  & 8323.2  & 4024  & 3353  & 3423  & 0.0   & 9     & 21    & \multicolumn{1}{r}{8392} & 107658 & 44936 & 39691 & \multicolumn{1}{r}{3.5 } & 3600  & 9273  & \multicolumn{1}{r}{1.01 } \\
    6     & 29    & 281   & 8278  & 8339.0  & 5080  & 3272  & 3351  & 0.0   & 15    & 45    & \multicolumn{1}{r}{8315} & 66442 & 48836 & 39775 & \multicolumn{1}{r}{3.1 } & 3600  & 9412  & \multicolumn{1}{r}{1.00 } \\
    6     & 32    & 281   & 8294  & 8358.2  & 13321  & 3340  & 3428  & 0.0   & 19    & 74    & \multicolumn{1}{r}{9778} & 101634 & 53076 & 39868 & \multicolumn{1}{r}{18.5 } & 3600  & 9637  & \multicolumn{1}{r}{1.18 } \\
    6     & 36    & 281   & 8449  & 8563.0  & 6405  & 3428  & 3528  & 0.0   & 20    & 93    & \multicolumn{1}{r}{10681} & 71224 & 59825 & 40006 & \multicolumn{1}{r}{25.6 } & 3600  & 9701  & \multicolumn{1}{r}{1.26 } \\
    6     & 41    & 281   & 8518  & 8720.2  & 18609  & 9476  & 9591  & 0.0   & 42    & 74    & \multicolumn{1}{r}{8260} & 12550 & 70639 & 40201 & \multicolumn{1}{r}{3.4 } & 3600  & 10144 & \multicolumn{1}{r}{0.97 } \\
    6     & 47    & 281   & 8739  & 8856.0  & 3942  & 9634  & 9767  & 0.0   & 24    & 79    & \multicolumn{1}{r}{8306} & 2972  & 87653 & 40468 & \multicolumn{1}{r}{4.5 } & 3600  & 10343 & \multicolumn{1}{r}{0.95 } \\
    6     & 57    & 281   & 9006  & 9090.0  & 352   & 9914  & 10077  & 0.0   & 4     & 63    & -     & 983   & 127690 & 40993 & -     & 3600  & 10865 & - \\
    \\
    7     & 25    & 361   & 9369  & 9461.6  & 5705  & 2387  & 2431  & 0.0   & 3     & 20    & \multicolumn{1}{r}{9504} & 78862 & 67532 & 65305 & \multicolumn{1}{r}{4.7 } & 3600  & 10009 & \multicolumn{1}{r}{1.01 } \\
    7     & 26    & 361   & 9536  & 9536.0  & 4284  & 2470  & 2516  & 0.0   & 3     & 8     & \multicolumn{1}{r}{9461} & 67573 & 68909 & 65331 & \multicolumn{1}{r}{4.1 } & 3600  & 10018 & \multicolumn{1}{r}{0.99 } \\
    7     & 28    & 361   & 9448  & 9521.2  & 6659  & 2460  & 2510  & 0.0   & 12    & 23    & \multicolumn{1}{r}{9579} & 67696 & 71166 & 65386 & \multicolumn{1}{r}{5.1 } & 3600  & 10050 & \multicolumn{1}{r}{1.01 } \\
    7     & 31    & 361   & 9499  & 9586.0  & 7977  & 3428  & 3513  & 0.0   & 14    & 43    & \multicolumn{1}{r}{9604} & 39028 & 75396 & 65476 & \multicolumn{1}{r}{4.2 } & 3600  & 10840 & \multicolumn{1}{r}{1.01 } \\
    7     & 33    & 361   & 9507  & 9629.0  & 4903  & 3334  & 3425  & 0.0   & 19    & 59    & \multicolumn{1}{r}{9467} & 13292 & 78675 & 65541 & \multicolumn{1}{r}{3.2 } & 3600  & 10956 & \multicolumn{1}{r}{1.00 } \\
    7     & 37    & 361   & 9717  & 9762.8  & 6950  & 3372  & 3475  & 0.0   & 26    & 63    & -     & 30574 & 85902 & 65683 & -     & 3600  & 11085 & - \\
    7     & 41    & 361   & 9720  & 9746.0  & 4653  & 9397  & 9512  & 0.0   & 16    & 61    & -     & 16766 & 94861 & 65841 & -     & 3600  & 11100 & - \\
    7     & 46    & 361   & 9832  & 9891.0  & 11250  & 9519  & 9649  & 0.0   & 24    & 77    & -     & 6419  & 109364 & 66061 & -     & 3600  & 11367 & - \\
    7     & 52    & 361   & 9793  & 9918.8  & 4272  & 9585  & 9733  & 0.0   & 17    & 131   & -     & 3155  & 131150 & 66358 & -     & 3600  & 11465 & - \\
    7     & 61    & 361   & 9998  & 10011.0  & 5210  & 9802  & 9977  & 0.0   & 30    & 324   & -     & 497   & 173433 & 66871 & -     & 3600  & 11909 & - \\
    7     & 73    & 361   & 9911  & 10076.2  & 1475  & 9779  & 9990  & 0.0   & 4     & 143   & -     & 1     & 252900 & 67681 & -     & 3600  & 12431 & - \\
    \\
    9     & 18    & 256   & 281   & 282.4  & 4863  & 2991  & 3037  & 0.0   & 2     & 5     & \multicolumn{1}{r}{281} & 170032 & 33099 & 32811 & \multicolumn{1}{r}{0.0 } & 2061  & 307   & \multicolumn{1}{r}{1.00 } \\
    9     & 19    & 256   & 283   & 283.0  & 3449  & 3419  & 3468  & 0.0   & 5     & 16    & \multicolumn{1}{r}{279} & 51107 & 33688 & 32830 & \multicolumn{1}{r}{0.0 } & 651   & 307   & \multicolumn{1}{r}{0.99 } \\
    9     & 20    & 256   & 280   & 280.0  & 1988  & 3440  & 3492  & 0.0   & 6     & 20    & -     & 288725 & 34344 & 32850 & -     & 3600  & 302   & - \\
    9     & 22    & 256   & 276   & 278.0  & 8725  & 3372  & 3430  & 0.0   & 14    & 23    & \multicolumn{1}{r}{276} & 409845 & 35761 & 32893 & \multicolumn{1}{r}{0.4 } & 3600  & 303   & \multicolumn{1}{r}{1.00 } \\
    9     & 24    & 256   & 276   & 279.8  & 3969  & 3384  & 3448  & 0.0   & 7     & 28    & \multicolumn{1}{r}{276} & 450126 & 37419 & 32940 & \multicolumn{1}{r}{0.0 } & 1270  & 300   & \multicolumn{1}{r}{1.00 } \\
    9     & 26    & 256   & 278   & 279.0  & 9238  & 3375  & 3445  & 0.0   & 22    & 62    & \multicolumn{1}{r}{339} & 165847 & 39431 & 32991 & \multicolumn{1}{r}{20.1 } & 3600  & 311   & \multicolumn{1}{r}{1.22 } \\
    9     & 29    & 256   & 279   & 282.0  & 7038  & 3330  & 3409  & 0.0   & 22    & 189   & \multicolumn{1}{r}{273} & 64689 & 42851 & 33075 & \multicolumn{1}{r}{0.0 } & 2436  & 321   & \multicolumn{1}{r}{0.98 } \\
    9     & 32    & 256   & 284   & 286.0  & 6359  & 3453  & 3541  & 0.0   & 21    & 164   & \multicolumn{1}{r}{273} & 83555 & 46921 & 33168 & \multicolumn{1}{r}{0.0 } & 2063  & 321   & \multicolumn{1}{r}{0.96 } \\
    9     & 37    & 256   & 283   & 289.0  & 3500  & 3199  & 3302  & 0.0   & 18    & 85    & \multicolumn{1}{r}{343} & 88225 & 55610 & 33343 & \multicolumn{1}{r}{20.4 } & 3600  & 335   & \multicolumn{1}{r}{1.21 } \\
    9     & 43    & 256   & 286   & 289.2  & 2051  & 9563  & 9684  & 0.0   & 9     & 47    & -     & 73001 & 69711 & 33586 & -     & 3600  & 352   & - \\
    9     & 52    & 256   & 284   & 288.0  & 2943  & 9776  & 9850  & 0.0   & 18    & 117   & -     & 14922 & 99686 & 34018 & -     & 3600  & 347   & - \\
    \bottomrule
    \end{tabular}}%
}
{}
\end{table}%

\begin{table}[htbp]
  \TABLE
  {{Detailed results for the SoftCluVRP instances (Part2) }  \label{tab: SoftCluE2}}
  {\resizebox{\columnwidth}{!}{%
    \begin{tabular}{lllrrrrrrrrrrrrrrrr}
        \toprule
        \multicolumn{1}{c}{\multirow{2}[0]{*}{$\Delta$}}  &\multicolumn{1}{c}{\multirow{2}[0]{*}{Cluster}} &\multicolumn{1}{c}{\multirow{2}[0]{*}{Customer}}  & \multicolumn{8}{c}{SIM} & \multicolumn{6}{c}{Exact}& &  \\
        \cmidrule(lr){4-11}
        \cmidrule(lr){12-17}
        \cmidrule(lr){18-19}
        & & &$\text{Obj}_{\text{best}}$ &$\text{Obj}_{\text{avg}}$ &\multicolumn{1}{l}{Nodes} & \multicolumn{1}{l}{Rows} & \multicolumn{1}{l}{Cols} & \multicolumn{1}{l}{Gap($\%$)} & $\text{Time}_{\text{MIP}}$ & $\text{Time}_{\text{Total}}$& Obj&Nodes &Rows &Cols & Gap($\%$) &$\text{Time}_{\text{MIP}}$& \multicolumn{1}{l}{Min(sample)}& $\text{Exact/SIM}_{\text{best}}$ \\
     \hline
    13    & 17    & 253   & 535   & 535.0  & 2014  & 3353  & 3396  & 0.0   & 2     & 12    & \multicolumn{1}{r}{530} & 51520 & 32049 & 32031 & \multicolumn{1}{r}{0.0 } & 407   & 548   & \multicolumn{1}{r}{0.99 } \\
    13    & 19    & 253   & 530   & 533.6  & 5863  & 3210  & 3259  & 0.0   & 5     & 36    & \multicolumn{1}{r}{521} & 7346  & 33220 & 32068 & \multicolumn{1}{r}{0.0 } & 207   & 543   & \multicolumn{1}{r}{0.98 } \\
    13    & 20    & 253   & 532   & 532.0  & 4597  & 3258  & 3310  & 0.0   & 3     & 37    & \multicolumn{1}{r}{521} & 51999 & 33820 & 32088 & \multicolumn{1}{r}{0.0 } & 489   & 549   & \multicolumn{1}{r}{0.98 } \\
    13    & 22    & 253   & 533   & 533.8  & 2319  & 3409  & 3467  & 0.0   & 4     & 16    & \multicolumn{1}{r}{523} & 50875 & 35209 & 32131 & \multicolumn{1}{r}{0.0 } & 840   & 558   & \multicolumn{1}{r}{0.98 } \\
    13    & 23    & 253   & 533   & 533.0  & 3292  & 3458  & 3519  & 0.0   & 4     & 24    & \multicolumn{1}{r}{523} & 63188 & 35965 & 32154 & \multicolumn{1}{r}{0.0 } & 675   & 560   & \multicolumn{1}{r}{0.98 } \\
    13    & 26    & 253   & 531   & 532.2  & 4990  & 3354  & 3424  & 0.0   & 7     & 27    & \multicolumn{1}{r}{523} & 132609 & 38705 & 32229 & \multicolumn{1}{r}{0.0 } & 1513  & 579   & \multicolumn{1}{r}{0.98 } \\
    13    & 29    & 253   & 534   & 537.0  & 6656  & 3316  & 3395  & 0.0   & 9     & 28    & \multicolumn{1}{r}{522} & 103643 & 42111 & 32313 & \multicolumn{1}{r}{0.0 } & 1584  & 580   & \multicolumn{1}{r}{0.98 } \\
    13    & 32    & 253   & 530   & 533.2  & 2363  & 3603  & 3691  & 0.0   & 4     & 17    & \multicolumn{1}{r}{521} & 84100 & 46222 & 32406 & \multicolumn{1}{r}{0.0 } & 1922  & 579   & \multicolumn{1}{r}{0.98 } \\
    13    & 37    & 253   & 533   & 537.8  & 3490  & 3369  & 3472  & 0.0   & 18    & 41    & \multicolumn{1}{r}{551} & 152979 & 55002 & 32581 & \multicolumn{1}{r}{5.6 } & 3600  & 593   & \multicolumn{1}{r}{1.03 } \\
    13    & 43    & 253   & 533   & 536.0  & 3930  & 9538  & 9659  & 0.0   & 11    & 84    & \multicolumn{1}{r}{521} & 67687 & 69112 & 32824 & \multicolumn{1}{r}{0.0 } & 2790  & 611   & \multicolumn{1}{r}{0.98 } \\
    13    & 51    & 253   & 538   & 547.0  & 4103  & 9463  & 9608  & 0.0   & 24    & 122   & \multicolumn{1}{r}{522} & 41502 & 95097 & 33204 & \multicolumn{1}{r}{0.4 } & 3600  & 631   & \multicolumn{1}{r}{0.97 } \\
    \\
    17    & 17    & 241   & 393   & 393.2  & 969   & 2427  & 2455  & 0.0   & 1     & 9     & \multicolumn{1}{r}{386} & 42444 & 28814 & 29073 & \multicolumn{1}{r}{0.0 } & 223   & 410   & \multicolumn{1}{r}{0.98 } \\
    17    & 18    & 241   & 393   & 393.0  & 934   & 2342  & 2372  & 0.0   & 1     & 8     & \multicolumn{1}{r}{385} & 20643 & 29354 & 29091 & \multicolumn{1}{r}{0.0 } & 235   & 407   & \multicolumn{1}{r}{0.98 } \\
    17    & 19    & 241   & 388   & 388.0  & 2866  & 2358  & 2390  & 0.0   & 2     & 10    & \multicolumn{1}{r}{385} & 27296 & 29962 & 29110 & \multicolumn{1}{r}{0.0 } & 305   & 415   & \multicolumn{1}{r}{0.99 } \\
    17    & 21    & 241   & 388   & 388.0  & 4108  & 2360  & 2396  & 0.0   & 5     & 15    & \multicolumn{1}{r}{385} & 38918 & 31237 & 29151 & \multicolumn{1}{r}{0.0 } & 316   & 421   & \multicolumn{1}{r}{0.99 } \\
    17    & 22    & 241   & 388   & 388.0  & 4389  & 2364  & 2402  & 0.0   & 9     & 23    & \multicolumn{1}{r}{385} & 17394 & 31926 & 29173 & \multicolumn{1}{r}{0.0 } & 267   & 425   & \multicolumn{1}{r}{0.99 } \\
    17    & 25    & 241   & 387   & 389.4  & 8867  & 2414  & 2458  & 0.0   & 11    & 41    & \multicolumn{1}{r}{382} & 21980 & 34680 & 29245 & \multicolumn{1}{r}{0.0 } & 356   & 434   & \multicolumn{1}{r}{0.99 } \\
    17    & 27    & 241   & 385   & 389.0  & 2483  & 2570  & 2618  & 0.0   & 4     & 19    & \multicolumn{1}{r}{382} & 77815 & 36786 & 29298 & \multicolumn{1}{r}{0.0 } & 1036  & 416   & \multicolumn{1}{r}{0.99 } \\
    17    & 31    & 241   & 392   & 393.0  & 4195  & 3319  & 3404  & 0.0   & 13    & 28    & \multicolumn{1}{r}{445} & 232496 & 41893 & 29416 & \multicolumn{1}{r}{13.3 } & 3600  & 470   & \multicolumn{1}{r}{1.14 } \\
    17    & 35    & 241   & 397   & 402.0  & 32335  & 3531  & 3628  & 0.0   & 143   & 408   & \multicolumn{1}{r}{390} & 62163 & 48367 & 29550 & \multicolumn{1}{r}{0.3 } & 3600  & 471   & \multicolumn{1}{r}{0.98 } \\
    17    & 41    & 241   & 394   & 397.0  & 10248  & 9165  & 9280  & 0.0   & 132   & 294   & \multicolumn{1}{r}{607} & 65545 & 61041 & 29781 & \multicolumn{1}{r}{36.6 } & 3600  & 490   & \multicolumn{1}{r}{1.54 } \\
    17    & 49    & 241   & 394   & 398.0  & 5129  & 9574  & 9713  & 0.0   & 20    & 94    & \multicolumn{1}{r}{524} & 23990 & 84779 & 30145 & \multicolumn{1}{r}{26.7 } & 3600  & 510   & \multicolumn{1}{r}{1.33 } \\
    \\
    18    & 21    & 301   & 566   & 566.0  & 7527  & 3337  & 3392  & 0.0   & 11    & 20    & \multicolumn{1}{r}{558} & 5435  & 46087 & 45381 & \multicolumn{1}{r}{0.0 } & 325   & 613   & \multicolumn{1}{r}{0.99 } \\
    18    & 22    & 301   & 566   & 566.0  & 1199  & 3606  & 3664  & 0.0   & 2     & 11    & \multicolumn{1}{r}{558} & 12606 & 46824 & 45403 & \multicolumn{1}{r}{0.0 } & 463   & 624   & \multicolumn{1}{r}{0.99 } \\
    18    & 24    & 301   & 566   & 566.0  & 4220  & 3383  & 3447  & 0.0   & 10    & 23    & \multicolumn{1}{r}{558} & 10086 & 48415 & 45450 & \multicolumn{1}{r}{0.0 } & 424   & 639   & \multicolumn{1}{r}{0.99 } \\
    18    & 26    & 301   & 566   & 569.0  & 19255  & 3792  & 3864  & 0.0   & 61    & 86    & \multicolumn{1}{r}{562} & 81781 & 50828 & 45501 & \multicolumn{1}{r}{0.0 } & 1922  & 642   & \multicolumn{1}{r}{0.99 } \\
    18    & 28    & 301   & 567   & 570.0  & 6717  & 3349  & 3425  & 0.0   & 17    & 74    & \multicolumn{1}{r}{621} & 193951 & 53317 & 45556 & \multicolumn{1}{r}{11.6 } & 3600  & 629   & \multicolumn{1}{r}{1.10 } \\
    18    & 31    & 301   & 566   & 566.0  & 8137  & 3373  & 3458  & 0.0   & 22    & 213   & \multicolumn{1}{r}{554} & 81956 & 57307 & 45646 & \multicolumn{1}{r}{0.0 } & 3233  & 644   & \multicolumn{1}{r}{0.98 } \\
    18    & 34    & 301   & 573   & 576.6  & 5543  & 3588  & 3682  & 0.0   & 12    & 28    & \multicolumn{1}{r}{568} & 83856 & 62106 & 45745 & \multicolumn{1}{r}{3.2 } & 3600  & 636   & \multicolumn{1}{r}{0.99 } \\
    18    & 38    & 301   & 569   & 574.0  & 1867  & 3659  & 3765  & 0.0   & 7     & 29    & \multicolumn{1}{r}{764} & 93716 & 69755 & 45891 & \multicolumn{1}{r}{28.3 } & 3600  & 670   & \multicolumn{1}{r}{1.34 } \\
    18    & 43    & 301   & 576   & 581.8  & 3258  & 9607  & 9728  & 0.0   & 13    & 110   & \multicolumn{1}{r}{563} & 10279 & 81973 & 46096 & \multicolumn{1}{r}{2.5 } & 3600  & 678   & \multicolumn{1}{r}{0.98 } \\
    18    & 51    & 301   & 566   & 588.0  & 10806  & 9576  & 9721  & 0.0   & 64    & 193   & \multicolumn{1}{r}{566} & 6558  & 108039 & 46476 & \multicolumn{1}{r}{3.4 } & 3600  & 707   & \multicolumn{1}{r}{1.00 } \\
    18    & 61    & 301   & 587   & 600.0  & 14139  & 9548  & 9723  & 0.0   & 93    & 414   & -     & 64    & 154177 & 47041 & -     & 3600  & 739   & - \\
    \\
    19    & 25    & 361   & 899   & 908.0  & 49785  & 9525  & 9730  & 26.5  & 1889  & 3600  & \multicolumn{1}{r}{888} & 147472 & 67339 & 65305 & \multicolumn{1}{r}{1.8 } & 3600  & 1012  & \multicolumn{1}{r}{0.99 } \\
    19    & 26    & 361   & 934   & 934.0  & 50623  & 8476  & 8690  & 13.1  & 3593  & 3600  & \multicolumn{1}{r}{894} & 153976 & 68552 & 65331 & \multicolumn{1}{r}{2.7 } & 3600  & 1041  & \multicolumn{1}{r}{0.96 } \\
    19    & 28    & 361   & 765   & 767.0  & 12458  & 3221  & 3297  & 0.0   & 35    & 127   & \multicolumn{1}{r}{741} & 76185 & 71214 & 65386 & \multicolumn{1}{r}{0.0 } & 3070  & 860   & \multicolumn{1}{r}{0.97 } \\
    19    & 31    & 361   & 766   & 768.8  & 2781  & 3516  & 3601  & 0.0   & 7     & 24    & \multicolumn{1}{r}{916} & 90675 & 75856 & 65476 & \multicolumn{1}{r}{21.0 } & 3600  & 873   & \multicolumn{1}{r}{1.20 } \\
    19    & 33    & 361   & 739   & 744.0  & 3071  & 3760  & 3851  & 0.0   & 5     & 129   & \multicolumn{1}{r}{727} & 9636  & 79362 & 65541 & \multicolumn{1}{r}{0.0 } & 3464  & 863   & \multicolumn{1}{r}{0.98 } \\
    19    & 37    & 361   & 746   & 767.2  & 12686  & 4518  & 4656  & 0.0   & 70    & 260   & \multicolumn{1}{r}{798} & 18486 & 86959 & 65683 & \multicolumn{1}{r}{9.3 } & 3600  & 893   & \multicolumn{1}{r}{1.07 } \\
    19    & 41    & 361   & 769   & 773.4  & 7469  & 9557  & 9711  & 0.0   & 37    & 204   & \multicolumn{1}{r}{927} & 16062 & 96064 & 65841 & \multicolumn{1}{r}{22.0 } & 3600  & 932   & \multicolumn{1}{r}{1.21 } \\
    19    & 46    & 361   & 769   & 777.0  & 4524  & 9720  & 9894  & 0.0   & 27    & 203   & \multicolumn{1}{r}{735} & 6480  & 109938 & 66061 & \multicolumn{1}{r}{1.6 } & 3600  & 929   & \multicolumn{1}{r}{0.96 } \\
    19    & 52    & 361   & 780   & 791.0  & 13053  & 9633  & 9831  & 0.0   & 114   & 256   & -     & 6370  & 131081 & 66358 & -     & 3600  & 1014  & - \\
    19    & 61    & 361   & 774   & 785.4  & 9080  & 9419  & 9653  & 0.0   & 144   & 453   & -     & 120   & 173470 & 66871 & -     & 3600  & 1035  & - \\
    19    & 73    & 361   & 825   & 835.0  & 8746  & 9652  & 9934  & 0.0   & 120   & 228   & -     & 1     & 253098 & 67681 & -     & 3600  & 1050  & - \\
    \hline  
     \multicolumn{3}{c}{Average(Table \ref{tab: SoftCluE1} \& \ref{tab: SoftCluE2})} & & & 7084.1 &5010.7&5100.3& & 72.0&146.7& &79438.6&63824.4&41447.3& (35/110) &\multicolumn{1}{r}{2881.0}& & 1.04 (64/110)\\
    \bottomrule
    \end{tabular}}%
}
{}
\end{table}%

\begin{table}[htbp]
  \TABLE
  {{The SIM in SoftCluVRP BKS and heuristic comparison (Part1)}  \label{tab: SoftCluH1}}
  {\resizebox{0.9\columnwidth}{!}{%
    \begin{tabular}{lllrrrlrrrrrr}
    \toprule
        \multicolumn{1}{c}{\multirow{2}[0]{*}{$\Delta$}}& \multicolumn{1}{c}{\multirow{2}[0]{*}{Cluster}}&  \multicolumn{1}{c}{\multirow{2}[0]{*}{Customer}} & \multicolumn{3}{c}{SIM(20s)} &  \multicolumn{1}{c}{\multirow{2}[0]{*}{BKS $\dag$}}&\multicolumn{3}{c}{LMNS}&\multicolumn{3}{c}{Ratio}   \\
        \cmidrule(lr){4-6}
        \cmidrule(lr){8-10}
        \cmidrule(lr){11-13}
        & & &$\text{Obj}_{\text{best}}$ &$\text{Obj}_{\text{avg}}$&$\text{Time}_{\text{Total}}$&  & $\text{Obj}_{\text{best}}$ &$\text{Obj}_{\text{avg}}$&$\text{Time}_{\text{Total}}$ & $\text{BKS/SIM}_{\text{best}}$ & $(\text{LMNS/SIM})_{\text{best}}$ & $(\text{LMNS/SIM})_{\text{avg}}$\\
        \hline
    1     & 17    & 241   & 4640  & 4640.0  & 6     & $4640^*$  & 4640  & 4640.0  & 30    & 1.00  & 1.00  & 1.00  \\
    1     & 18    & 241   & 4652  & 4664.6  & 9     & 4645  & 4645  & 4645.0  & 31    & 1.00  & 1.00  & 1.00  \\
    1     & 19    & 241   & 4667  & 4685.0  & 10    & 4650  & 4650  & 4650.0  & 33    & 1.00  & 1.00  & 0.99  \\
    1     & 21    & 241   & 4669  & 4681.0  & 13    & 4650  & 4650  & 4650.0  & 33    & 1.00  & 1.00  & 0.99  \\
    1     & 22    & 241   & 4677  & 4694.0  & 25    & 4650  & 4650  & 4650.0  & 33    & 0.99  & 0.99  & 0.99  \\
    1     & 25    & 241   & 4660  & 4664.6  & 30    & 4650  & 4650  & 4651.2 & 35    & 1.00  & 1.00  & 1.00  \\
    1     & 27    & 241   & 4658  & 4663.4  & 23    & 4652  & 4652  & 4652.0  & 35    & 1.00  & 1.00  & 1.00  \\
    1     & 31    & 241   & 4710  & 4716.0  & 49    & 4665  & 4665  & 4665.0  & 44    & 0.99  & 0.99  & 0.99  \\
    1     & 35    & 241   & 4732  & 4740.2  & 58    & 4619  & 4619  & 4619.8 & 46    & 0.98  & 0.98  & 0.97  \\
    1     & 41    & 241   & 4685  & 4726.0  & 121   & 4619  & 4619  & 4621.3 & 44    & 0.99  & 0.99  & 0.98  \\
    1     & 49    & 241   & 4803  & 4860.4  & 117   & 4607  & 4619  & 4625.5 & 47    & 0.96  & 0.96  & 0.95  \\
    \\
    2     & 22    & 321   & 7442  & 7485.6  & 24    & 7394  & 7394  & 7395.9 & 66    & 0.99  & 0.99  & 0.99  \\
    2     & 23    & 321   & 7392  & 7403.2  & 46    & $7369^*$  & 7372  & 7381.2 & 66    & 1.00  & 1.00  & 1.00  \\
    2     & 25    & 321   & 7414  & 7426.2  & 24    & 7367  & 7367  & 7370.4 & 69    & 0.99  & 0.99  & 0.99  \\
    2     & 27    & 321   & 7399  & 7402.4  & 25    & 7333  & 7333  & 7334.3 & 72    & 0.99  & 0.99  & 0.99  \\
    2     & 30    & 321   & 7428  & 7445.0  & 36    & 7329  & 7329  & 7329.0  & 78    & 0.99  & 0.99  & 0.98  \\
    2     & 33    & 321   & 7443  & 7493.8  & 48    & 7311  & 7311  & 7314.1 & 80    & 0.98  & 0.98  & 0.98  \\
    2     & 36    & 321   & 7396  & 7494.2  & 38    & 7293  & 7293  & 7293.2 & 84    & 0.99  & 0.99  & 0.97  \\
    2     & 41    & 321   & 7385  & 7451.0  & 69    & 7283  & 7283  & 7286.2 & 88    & 0.99  & 0.99  & 0.98  \\
    2     & 46    & 321   & 7449  & 7541.8  & 71    & 7284  & 7284  & 7290.7 & 95    & 0.98  & 0.98  & 0.97  \\
    2     & 54    & 321   & 7560  & 7720.0  & 122   & 7274  & 7277  & 7278.7 & 101   & 0.96  & 0.96  & 0.94  \\
    2     & 65    & 321   & 7697  & 7780.0  & 115   & 7261  & 7264  & 7272.4 & 104   & 0.94  & 0.94  & 0.93  \\
    \\
    5     & 14    & 201   & 6970  & 6970.0  & 19    & $6970^*$  & 6970  & 6970.0  & 22    & 1.00  & 1.00  & 1.00  \\
    5     & 15    & 201   & 6933  & 6933.0  & 5     & $6742^*$  & 6742  & 6752.0  & 26    & 0.97  & 0.97  & 0.97  \\
    5     & 16    & 201   & 6895  & 6920.2  & 10    & $6742^*$  & 6742  & 6849.1 & 26    & 0.98  & 0.98  & 0.99  \\
    5     & 17    & 201   & 6862  & 6862.0  & 5     & $6862^*$  & 6862  & 6868.0  & 26    & 1.00  & 1.00  & 1.00  \\
    5     & 19    & 201   & 7030  & 7030.0  & 28    & $6874^*$  & 6874  & 6874.0  & 25    & 0.98  & 0.98  & 0.98  \\
    5     & 21    & 201   & 6881  & 6885.0  & 20    & $6816^*$  & 6816  & 6817.4 & 26    & 0.99  & 0.99  & 0.99  \\
    5     & 23    & 201   & 6887  & 6943.8  & 17    & $6750^*$  & 6750  & 6750.0  & 25    & 0.98  & 0.98  & 0.97  \\
    5     & 26    & 201   & 6880  & 6882.6  & 9     & $6704^*$  & 6704  & 6704.0  & 27    & 0.97  & 0.97  & 0.97  \\
    5     & 29    & 201   & 6879  & 6910.0  & 25    & $6704^*$  & 6704  & 6704.0  & 28    & 0.97  & 0.97  & 0.97  \\
    5     & 34    & 201   & 6833  & 7031.0  & 26    & 6684  & 6684  & 6692.4 & 29    & 0.98  & 0.98  & 0.95  \\
    5     & 41    & 201   & 7097  & 7128.2  & 125   & $6557^*$  & 6557  & 6578.2 & 32    & 0.92  & 0.92  & 0.92  \\
    \\
    6     & 19    & 281   & 8115  & 8122.2  & 14    & $8115^*$  & 8115  & 8115.3 & 54    & 1.00  & 1.00  & 1.00  \\
    6     & 21    & 281   & 8119  & 8134.4  & 7     & $8119^*$  & 8119  & 8125.5 & 52    & 1.00  & 1.00  & 1.00  \\
    6     & 22    & 281   & 8107  & 8107.0  & 15    & $8107^*$  & 8107  & 8113.7 & 52    & 1.00  & 1.00  & 1.00  \\
    6     & 24    & 281   & 8316  & 8351.0  & 8     & $8316^*$  & 8316  & 8318.8 & 52    & 1.00  & 1.00  & 1.00  \\
    6     & 26    & 281   & 8285  & 8305.0  & 25    & $8249^*$  & 8249  & 8256.4 & 54    & 1.00  & 1.00  & 0.99  \\
    6     & 29    & 281   & 8278  & 8320.8  & 91    & 8244  & 8244  & 8251.4 & 60    & 1.00  & 1.00  & 0.99  \\
    6     & 32    & 281   & 8486  & 8504.0  & 69    & 8179  & 8179  & 8197.3 & 59    & 0.96  & 0.96  & 0.96  \\
    6     & 36    & 281   & 8331  & 8470.6  & 52    & 8179  & 8179  & 8180.9 & 59    & 0.98  & 0.98  & 0.97  \\
    6     & 41    & 281   & 8595  & 8679.8  & 100   & 8204  & 8204  & 8206.5 & 66    & 0.95  & 0.95  & 0.95  \\
    6     & 47    & 281   & 8686  & 8770.0  & 104   & 8179  & 8179  & 8192.6 & 65    & 0.94  & 0.94  & 0.93  \\
    6     & 57    & 281   & 9059  & 9116.0  & 307   & 8204  & 8204  & 8205.6 & 75    & 0.91  & 0.91  & 0.90  \\
    \\
    7     & 25    & 361   & 9369  & 9461.6  & 20    & 9318  & 9318  & 9321.5 & 99    & 0.99  & 0.99  & 0.99  \\
    7     & 26    & 361   & 9536  & 9536.0  & 8     & 9295  & 9307  & 9314.1 & 101   & 0.97  & 0.98  & 0.98  \\
    7     & 28    & 361   & 9448  & 9521.2  & 17    & 9271  & 9272  & 9282.7 & 109   & 0.98  & 0.98  & 0.97  \\
    7     & 31    & 361   & 9667  & 9752.0  & 25    & 9418  & 9418  & 9442.6 & 101   & 0.97  & 0.97  & 0.97  \\
    7     & 33    & 361   & 9687  & 9731.0  & 37    & 9395  & 9401  & 9401.8 & 103   & 0.97  & 0.97  & 0.97  \\
    7     & 37    & 361   & 9764  & 9776.2  & 66    & 9395  & 9395  & 9403.7 & 104   & 0.96  & 0.96  & 0.96  \\
    7     & 41    & 361   & 9729  & 9747.2  & 106   & 9386  & 9386  & 9400.3 & 108   & 0.96  & 0.96  & 0.96  \\
    7     & 46    & 361   & 9685  & 9771.0  & 149   & 9368  & 9368  & 9376.7 & 102   & 0.97  & 0.97  & 0.96  \\
    7     & 52    & 361   & 9788  & 9890.0  & 184   & 9365  & 9365  & 9373.1 & 114   & 0.96  & 0.96  & 0.95  \\
    7     & 61    & 361   & 10140 & 10185.8  & 587   & 9316  & 9316  & 9343.6 & 128   & 0.92  & 0.92  & 0.92  \\
    7     & 73    & 361   & 9911  & 10076.2  & 143   & 9302  & 9302  & 9314.9 & 145   & 0.94  & 0.94  & 0.92  \\
    \\
    9     & 18    & 256   & 281   & 282.4  & 5     & $281^*$   & 281   & 281.0  & 39    & 1.00  & 1.00  & 1.00  \\
    9     & 19    & 256   & 283   & 283.0  & 16    & $279^*$   & 279   & 279.2 & 38    & 0.99  & 0.99  & 0.99  \\
    9     & 20    & 256   & 280   & 280.0  & 20    & $276^*$   & 276   & 276.6 & 40    & 0.99  & 0.99  & 0.99  \\
    9     & 22    & 256   & 276   & 278.0  & 18    & $276^*$   & 276   & 276.7 & 44    & 1.00  & 1.00  & 1.00  \\
    9     & 24    & 256   & 276   & 279.8  & 28    & $276^*$   & 276   & 276.9 & 44    & 1.00  & 1.00  & 0.99  \\
    9     & 26    & 256   & 278   & 279.6  & 35    & $273^*$   & 273   & 273.9 & 46    & 0.98  & 0.98  & 0.98  \\
    9     & 29    & 256   & 279   & 281.0  & 35    & $273^*$   & 273   & 273.6 & 45    & 0.98  & 0.98  & 0.97  \\
    9     & 32    & 256   & 284   & 286.4  & 57    & $273^*$   & 273   & 273.9 & 48    & 0.96  & 0.96  & 0.96  \\
    9     & 37    & 256   & 289   & 290.4  & 27    & 273   & 273   & 273.9 & 50    & 0.94  & 0.94  & 0.94  \\
    9     & 43    & 256   & 287   & 291.0  & 51    & 270   & 270   & 270.8 & 53    & 0.94  & 0.94  & 0.93  \\
    9     & 52    & 256   & 293   & 296.0  & 49    & 269   & 269   & 269.0  & 57    & 0.92  & 0.92  & 0.91  \\
    \bottomrule
    \end{tabular}}%
}
{}
\end{table}%

\begin{table}[htbp]
  \TABLE
  {{The SIM in SoftCluVRP BKS and heuristic comparison (Part2)}  \label{tab: SoftCluH2}}
  {\resizebox{0.9\columnwidth}{!}{%
    \begin{tabular}{lllrrrlrrrrrr}
    \toprule
        \multicolumn{1}{c}{\multirow{2}[0]{*}{$\Delta$}}& \multicolumn{1}{c}{\multirow{2}[0]{*}{Cluster}}&  \multicolumn{1}{c}{\multirow{2}[0]{*}{Customer}} & \multicolumn{3}{c}{SIM(20s)} &  \multicolumn{1}{c}{\multirow{2}[0]{*}{BKS$\dag$}}&\multicolumn{3}{c}{LMNS}&\multicolumn{3}{c}{Ratio}   \\
        \cmidrule(lr){4-6}
        \cmidrule(lr){8-10}
        \cmidrule(lr){11-13}
        & & &$\text{Obj}_{\text{best}}$ &$\text{Obj}_{\text{avg}}$&$\text{Time}_{\text{Total}}$&  & $\text{Obj}_{\text{best}}$ &$\text{Obj}_{\text{avg}}$&$\text{Time}_{\text{Total}}$ & $\text{BKS/SIM}_{\text{best}}$ & $(\text{LMNS/SIM})_{\text{best}}$ & $(\text{LMNS/SIM})_{\text{avg}}$\\
        \hline
      13    & 17    & 253   & 535   & 535.0  & 12    & $530^*$   & 530   & 530.4 & 40    & 0.99  & 0.99  & 0.99  \\
    13    & 19    & 253   & 530   & 533.6  & 36    & $521^*$   & 521   & 521.8 & 40    & 0.98  & 0.98  & 0.98  \\
    13    & 20    & 253   & 532   & 532.0  & 37    & $521^*$   & 521   & 521.5 & 42    & 0.98  & 0.98  & 0.98  \\
    13    & 22    & 253   & 533   & 533.8  & 16    & $523^*$   & 523   & 523.2 & 42    & 0.98  & 0.98  & 0.98  \\
    13    & 23    & 253   & 533   & 533.0  & 24    & $523^*$   & 523   & 523.2 & 43    & 0.98  & 0.98  & 0.98  \\
    13    & 26    & 253   & 531   & 532.2  & 27    & $523^*$   & 523   & 523.0  & 46    & 0.98  & 0.98  & 0.98  \\
    13    & 29    & 253   & 534   & 537.0  & 28    & $522^*$   & 522   & 522.0  & 48    & 0.98  & 0.98  & 0.97  \\
    13    & 32    & 253   & 530   & 533.0  & 17    & $521^*$   & 521   & 521.2 & 49    & 0.98  & 0.98  & 0.98  \\
    13    & 37    & 253   & 533   & 537.8  & 41    & $521^*$   & 521   & 521.9 & 53    & 0.98  & 0.98  & 0.97  \\
    13    & 43    & 253   & 533   & 536.0  & 84    & 521   & 521   & 521.0  & 54    & 0.98  & 0.98  & 0.97  \\
    13    & 51    & 253   & 547   & 551.0  & 78    & 521   & 521   & 521.0  & 58    & 0.95  & 0.95  & 0.95  \\
    \\
    17    & 17    & 241   & 393   & 393.2  & 9     & $386^*$   & 386   & 386.0  & 44    & 0.98  & 0.98  & 0.98  \\
    17    & 18    & 241   & 393   & 393.0  & 8     & $385^*$   & 385   & 385.0  & 45    & 0.98  & 0.98  & 0.98  \\
    17    & 19    & 241   & 388   & 388.0  & 10    & $385^*$   & 385   & 385.0  & 46    & 0.99  & 0.99  & 0.99  \\
    17    & 21    & 241   & 388   & 388.0  & 15    & $385^*$   & 385   & 385.0  & 47    & 0.99  & 0.99  & 0.99  \\
    17    & 22    & 241   & 388   & 388.0  & 23    & $385^*$   & 385   & 385.0  & 47    & 0.99  & 0.99  & 0.99  \\
    17    & 25    & 241   & 387   & 389.4  & 41    & $382^*$   & 382   & 382.2 & 47    & 0.99  & 0.99  & 0.98  \\
    17    & 27    & 241   & 385   & 389.0  & 19    & $382^*$   & 382   & 382.0  & 49    & 0.99  & 0.99  & 0.98  \\
    17    & 31    & 241   & 392   & 394.8  & 24    & $390^*$   & 390   & 390.0  & 51    & 0.99  & 0.99  & 0.99  \\
    17    & 35    & 241   & 394   & 403.0  & 26    & 390   & 390   & 390.0  & 57    & 0.99  & 0.99  & 0.97  \\
    17    & 41    & 241   & 398   & 399.2  & 36    & 388   & 388   & 388.4 & 59    & 0.97  & 0.97  & 0.97  \\
    17    & 49    & 241   & 414   & 416.4  & 70    & 387   & 387   & 387.2 & 60    & 0.93  & 0.93  & 0.93  \\
    \\
    18    & 21    & 301   & 566   & 566.0  & 16    & $558^*$   & 558   & 558.0  & 58    & 0.99  & 0.99  & 0.99  \\
    18    & 22    & 301   & 566   & 566.0  & 11    & $558^*$   & 558   & 558.0  & 59    & 0.99  & 0.99  & 0.99  \\
    18    & 24    & 301   & 566   & 566.0  & 46    & $558^*$   & 558   & 558.0  & 64    & 0.99  & 0.99  & 0.99  \\
    18    & 26    & 301   & 573   & 574.8  & 58    & $562^*$   & 562   & 562.0  & 63    & 0.98  & 0.98  & 0.98  \\
    18    & 28    & 301   & 567   & 572.0  & 81    & 558   & 558   & 558.0  & 66    & 0.98  & 0.98  & 0.98  \\
    18    & 31    & 301   & 566   & 566.0  & 36    & $554^*$   & 554   & 554.0  & 71    & 0.98  & 0.98  & 0.98  \\
    18    & 34    & 301   & 577   & 578.2  & 35    & $554^*$   & 554   & 554.1 & 70    & 0.96  & 0.96  & 0.96  \\
    18    & 38    & 301   & 579   & 597.0  & 40    & 555   & 555   & 555.1 & 74    & 0.96  & 0.96  & 0.93  \\
    18    & 43    & 301   & 588   & 591.8  & 36    & 558   & 558   & 558.0  & 83    & 0.95  & 0.95  & 0.94  \\
    18    & 51    & 301   & 619   & 620.0  & 54    & 555   & 555   & 555.9 & 83    & 0.90  & 0.90  & 0.90  \\
    18    & 61    & 301   & 622   & 628.0  & 75    & 556   & 556   & 556.6 & 92    & 0.89  & 0.89  & 0.89  \\
    \\
    19    & 25    & 361   & 928   & 931.6  & 25    & $886^*$   & 887   & 887.9 & 50    & 0.95  & 0.96  & 0.95  \\
    19    & 26    & 361   & 931   & 933.8  & 48    & $888^*$   & 889   & 889.0  & 51    & 0.95  & 0.95  & 0.95  \\
    19    & 28    & 361   & 763   & 766.0  & 47    & $741^*$   & 741   & 742.0  & 77    & 0.97  & 0.97  & 0.97  \\
    19    & 31    & 361   & 763   & 769.2  & 36    & 735   & 737   & 737.5 & 84    & 0.96  & 0.97  & 0.96  \\
    19    & 33    & 361   & 759   & 766.0  & 70    & $727^*$   & 728   & 729.1 & 89    & 0.96  & 0.96  & 0.95  \\
    19    & 37    & 361   & 751   & 772.6  & 58    & $732^*$   & 733   & 733.5 & 100   & 0.97  & 0.98  & 0.95  \\
    19    & 41    & 361   & 778   & 783.8  & 53    & 730   & 730   & 730.7 & 109   & 0.94  & 0.94  & 0.93  \\
    19    & 46    & 361   & 778   & 785.0  & 71    & 730   & 730   & 730.7 & 115   & 0.94  & 0.94  & 0.93  \\
    19    & 52    & 361   & 781   & 801.0  & 91    & 730   & 730   & 730.8 & 120   & 0.93  & 0.93  & 0.91  \\
    19    & 61    & 361   & 790   & 830.8  & 104   & 737   & 737   & 738.5 & 120   & 0.93  & 0.93  & 0.89  \\
    19    & 73    & 361   & 801   & 827.0  & 130   & 736   & 736   & 736.9 & 135   & 0.92  & 0.92  & 0.89  \\
    \hline
    \multicolumn{3}{c}{Average(Table \ref{tab: SoftCluH1} \& \ref{tab: SoftCluH2})} & & & 52.2 & & & &62.5& 0.97&0.97&0.97\\
    \bottomrule
    % \multicolumn{13}{l}{\smallskip
    % {\footnotesize $\dag$ The BKS objective values attached with $^*$ are proved to be optimal by \cite{hintsch2020exact}}}
    \end{tabular}}%
}
{\footnotesize $\dag$ The BKS objective values attached with $^*$ are proved to be optimal by \cite{hintsch2020exact}}
\end{table}%

% Table generated by Excel2LaTeX from sheet 'AP_result'

\subsubsection{Experiment results}
Table \ref{tab: SoftCluE1} and \ref{tab: SoftCluE2} present the details and solution outcomes for the SIM and the Exact. To verify the accuracy of our reproduced results, we record the number of instances solved to optimal. Among all considered instances, 35 achieved the optimal solution. This surpasses the original result of 21 by \cite{hessler2021branch}, indicating that our reproduction is not inferior to the original performance. We can measure the effectiveness of the SIM by the win rate, which denotes the proportion of instances in each group where the SIM's results surpassed that of the Exact. Out of 110 cases, SIM's best objective values outperformed the Exact in 64 instances. Moreover, Exact's best objective value is on average 4\% larger than that of the SIM. Finally, In comparison to the Exact under a 3600-second time limit, the SIM requires 146.7 seconds of computation time on average, which is only 5\% of that of the Exact. These results underscore the significant accuracy and efficiency of SIM in addressing this problem. Thus the questions \textbf{Q1} and \textbf{Q2} are answered.

\begin{figure}[htbp!]
	\centering
	\begin{minipage}{.35\textwidth}
		\FIGURE
		{\includegraphics[scale=0.35]{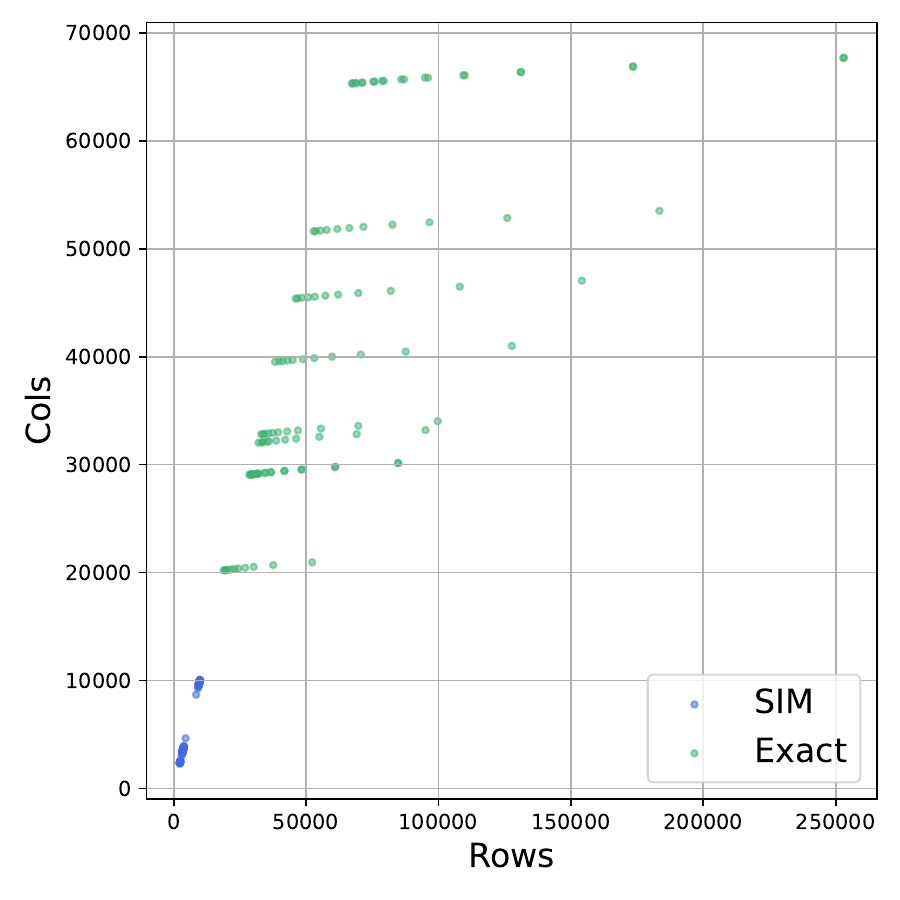}}
      %\captionsetup{width=0.9\textwidth}
		{{Scatter of model sizes}\label{fig: VRPscatter}}
        {}
	\end{minipage}%
	\begin{minipage}{.65\textwidth}
		\FIGURE
        {\includegraphics[scale=0.35]{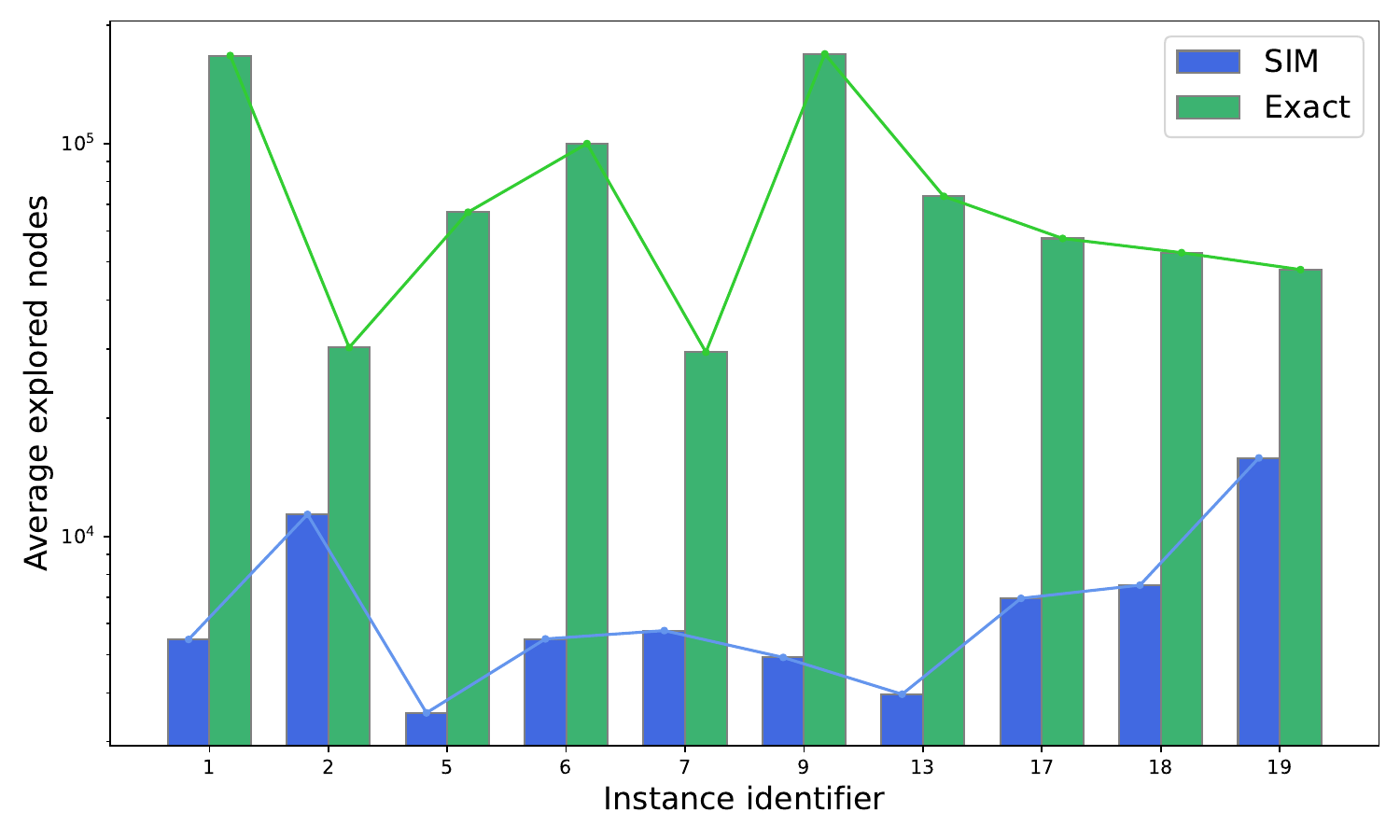}}
            %\captionsetup{width=0.9\textwidth}
		{{Explored nodes during the B\&C process}\label{fig: VRPnodes}}
        {}
	\end{minipage}
\end{figure}

Following the analysis approach used in Case I, we show the model details in Figure \ref{fig: VRPscatter}. It is evident from the figure that the model size of the SIM is concentrated within a smaller range of 2000 to 10000, whereas that of the Exact is distributed over a much larger interval, spanning from 20,000 to 250000. Regarding the average explored nodes within each instance group in Figure \ref{fig: VRPnodes}, due to the varying map of customers across different instance groups, the relationship between the explored nodes and the number of customers is not obverse. However, the figure clearly shows that the number of nodes in SIM is consistently and significantly lower than that in the Exact. On average, the explored nodes in SIM are only 9 \% of that in the Exact.

\begin{figure}[htbp!]
	\centering
	\begin{minipage}{.35\textwidth}
	\FIGURE
	{\includegraphics[scale=0.24]{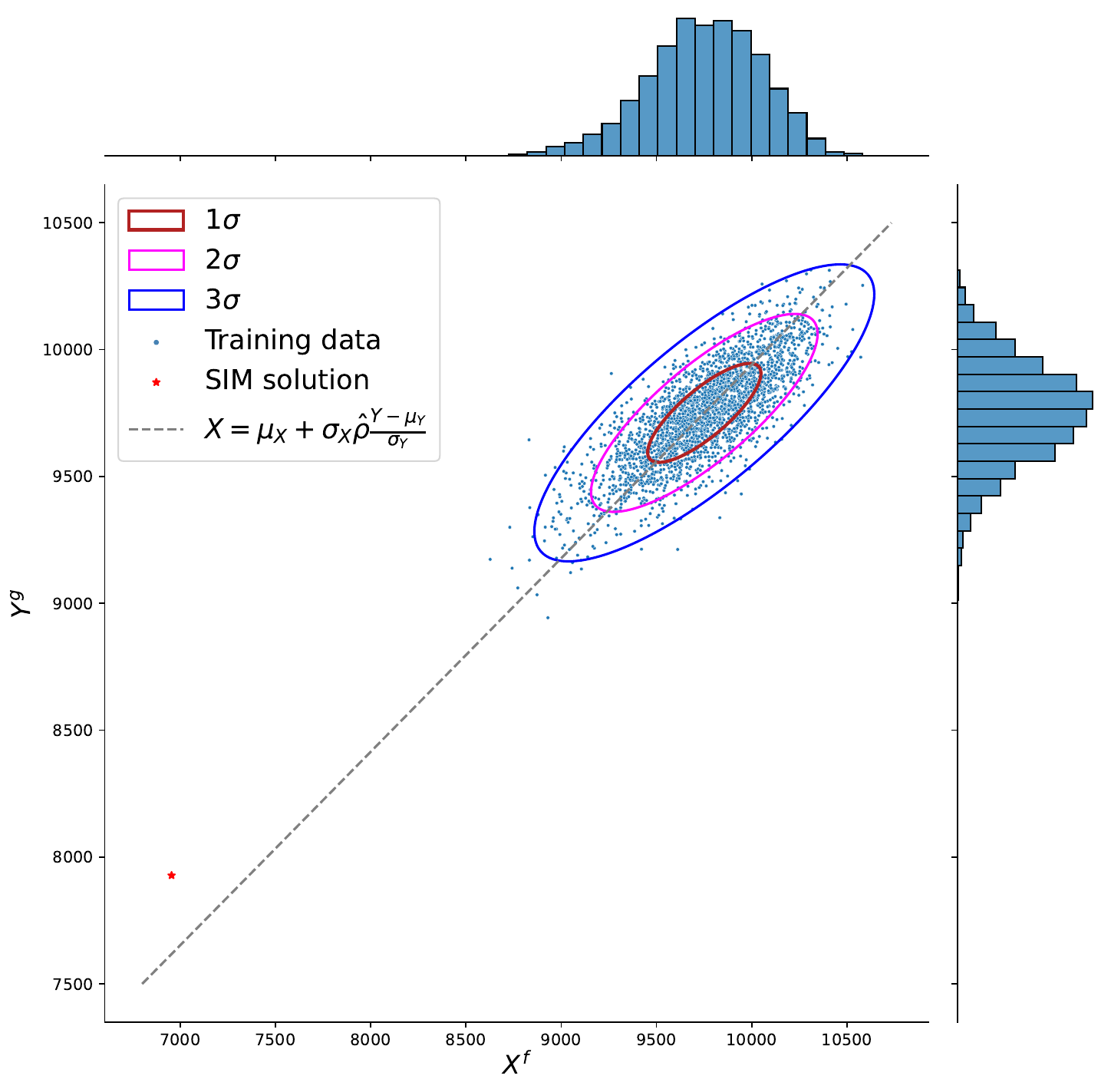}}
      %\captionsetup{width=0.9\textwidth}
    %\captionsetup{width=.9\linewidth}
    {{Parity plot of sampled data of the SIM and Exact models}\label{fig: VRPparity}}
    {}
	\end{minipage}%
	\begin{minipage}{.65\textwidth}
	\FIGURE
    {\includegraphics[scale=0.33]{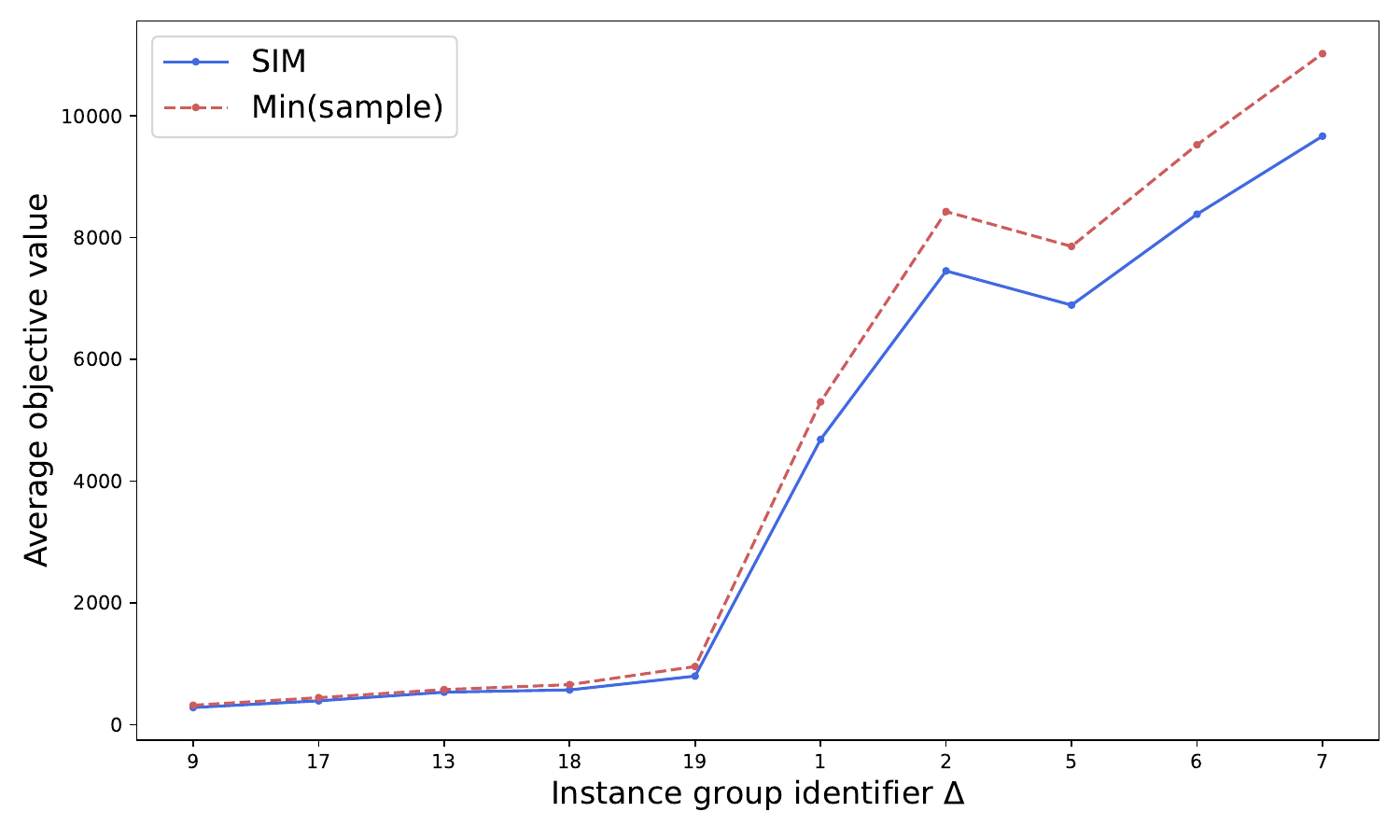}}
    %\captionsetup{width=0.9\textwidth}
    %\captionsetup{width=.9\linewidth}
	{{The comparison of average group objective values of the SIM the Min(sample)}\label{fig: VRPsample}}
    {}
	\end{minipage}
\end{figure}

In Figure \ref{fig: VRPparity}, we visualize the training data of a specific instance as a typical example. The training data in this figure exhibits a two-dimensional normal distribution. The optimal solution obtained by the SIM is located near the estimated line, suggesting that our theoretical approach applies to this problem. Figure \ref{fig: VRPsample} illustrates the comparison of the average of the minimum objective values obtained from different groups of case studies with the average objective values solved by the SIM. We can observe that the solutions derived from SIM consistently outperform the minimum values obtained from sampling. This also indicates that our framework exhibits certain learning capabilities in this problem.

The experimental data presented above demonstrates properties similar to those observed in Case I. Our framework is capable of generating models of limited size for this problem, which exhibit fewer explored nodes during the solving process. Furthermore, the outcomes display a notable difference with minimal sampling values, indicating a learning capacity.  It is suggested we can attribute the improvements of the SIM to the modifications by our framework, as an answer to \textbf{Q3}.

Finally, we restrict the computation time for each surrogate model within the Algorithm \ref{ag: train} to 20 seconds. This approach is then compared with the best-known solution (BKS) as documented in \cite{hintsch2020exact} and \cite{hintsch2021large}. Furthermore, a comparison is made with the latest state-of-the-art heuristic, referred to as LMNS in \cite{hintsch2021large}. The comparative outcomes are depicted in Tables \ref{tab: SoftCluH1} and \ref{tab: SoftCluH2}. We can observe that the objective values of the SIM are closely aligned with those of the BKS and LMNS, with all differences being within 3\% on average. Their results are almost equivalent for instances with no more than 30 clusters. Remarkably, the SIM demonstrated a superior average processing speed compared to the LMNS. This competitive performance suggests that the SIM could be considered as generally practical for the medium-sized SoftCluVRP as a heuristic algorithm.

\section{Conclusion}
\label{sec: con}
This study concentrates on the GAPR and presents a framework to generate surrogate models, improving the efficiency of B\&C-based problem-solving. To address the complexity involved by the routing constraints, we define a class of surrogate models to reduce decision variables and reconstruct constraints. We propose a framework for the learning of surrogate models through empirical data. Following this, the study provides theoretical insights concerning the representational power and statistical properties of the framework, yielding theoretical guarantees for its effectiveness. Experimental evaluations highlight the framework's accuracy and improved efficiency across varying scales within two practical problem classes, compared with the exact formulations as its learning target. Moreover, the framework shows close or superior performance compared with the state-of-the-art heuristics, providing supporting evidence on its general practicality, especially on small and medium problem sizes. Therefore, we can conclude that the proposed framework on the GAPR is effective in producing high-quality solutions with promoted efficiency.

In the future, the framework could be expanded to solve black-box optimization tasks. Furthermore, there is potential for integrating the framework with heuristic approaches or exact algorithms, serving as a warm start to facilitate the solution of complex problems.

\section*{Acknowledgments}
This work was funded by the National Natural Science Foundation of China under Grant No. 12320101001 and 12071428.

% \section*{Declaration of interest}
% No potential conflict of interest was reported by the authors

\bibliographystyle{informs2014}
\bibliography{references}
\end{document}